\documentclass[11pt,a4paper, final, twoside]{article}
\usepackage[style=numeric-comp, backend=bibtex, doi=false, isbn=false, natbib=true]{biblatex}
\usepackage{filecontents}
\usepackage[centertags]{amsmath}
\usepackage{amsfonts}
\usepackage{amssymb}
\usepackage{eqparbox}
\usepackage{amsthm}
\usepackage{newlfont}
\usepackage{fancyhdr}
\usepackage{amscd}
\usepackage{scalerel}
\usepackage{graphicx}
\usepackage{subfigure}
\usepackage{afterpage}
\usepackage{relsize}
\usepackage{sectsty}
\usepackage{siunitx}
\usepackage{etoolbox}
\usepackage{caption}
\usepackage{array}
\usepackage{booktabs}
\usepackage[table]{xcolor}
\usepackage{xcolor,colortbl}
\usepackage{colortbl,multirow,hhline}
\definecolor{blue(ncs)}{rgb}{0.0, 0.53, 0.74}
\usepackage[colorlinks=true, urlcolor=blue(ncs), linkcolor=blue(ncs), citecolor=blue(ncs)]{hyperref}
\usepackage{latexsym,epstopdf,float}
\usepackage{euscript,pifont,mathrsfs,latexsym,graphicx}
\usepackage{mathrsfs}
\usepackage{etoolbox}
\patchcmd\thebibliography{\labelsep}{\labelsep\itemsep=0pt\parsep=0pt\relax}{}{\typeout{Couldn't patch the command}}

\theoremstyle{plain}
\newtheorem{thm}{Theorem}[section]
\newtheorem{cor}[thm]{Corollary}
\newtheorem{lem}[thm]{Lemma}

\theoremstyle{definition}
\newtheorem{defn}{Definition}[section]
\newtheorem*{assum}{Assumptions}
\theoremstyle{remark}
\newtheorem{rem}{Remark}[section]
\newtheorem{exa}[rem]{Example}

 \numberwithin{equation}{section}
 \newcommand*{\medcup}{\mathbin{\scalebox{1.3}{\ensuremath{\bigcup}}}}
 \newcommand*{\sGamma}{\mathbin{\scalebox{.5}{\ensuremath{\Gamma}}}}

\begin{document}

\title{Optimality conditions for robust nonsmooth multiobjective optimization problems in Asplund spaces}
\author{Maryam Saadati, Morteza Oveisiha$^*$}
\date{}

\maketitle
\begin{center}
Department of Pure Mathematics, Faculty of Science, Imam Khomeini International University, P.O. Box 34149-16818, Qazvin, Iran.
\\
E-mail: m.saadati@edu.ikiu.ac.ir, oveisiha@sci.ikiu.ac.ir
\end{center}

\maketitle
\begin{abstract}
We employ a fuzzy optimality condition for the Fr\'{e}chet subdifferential and some advanced techniques of variational analysis such as formulae for the subdifferentials of an infinite family of nonsmooth functions and the coderivative scalarization to investigate robust optimality condition and robust duality for a nonsmooth/nonconvex multiobjective optimization problem dealing with uncertain constraints in arbitrary Asplund spaces. We establish necessary optimality conditions for weakly and properly robust efficient solutions of the problem in terms of the Mordukhovich subdifferentials of the related functions. Further, sufficient conditions for weakly and properly robust efficient solutions as well as for robust efficient solutions of the problem are provided by presenting new concepts of generalized convexity. Finally, we formulate a Mond-Weir-type robust dual problem to the reference problem, and examine weak, strong, and converse duality relations between them under the pseudo convexity assumptions.
\end{abstract}
\textbf{Keywords}\hspace{3mm}Robust nonsmooth multiobjective optimization . Optimality conditions . \mbox{Duality .} Limiting subdifferential . Generalized convexity
\newline
\textbf{Mathematics Subject Classification (2020)}\hspace{3mm}65K10 . 90C29 . 90C46

\footnote{\textsf{$^*$Corresponding author}}

\afterpage{} \fancyhead{} \fancyfoot{} \fancyhead[LE, RO]{\bf\thepage} \fancyhead[LO]{\small Optimality conditions for robust multiobjective optimization problems in Asplund spaces} \fancyhead[RE]{\small Saadati and Oveisiha}

\section{Introduction}\label{Sec1-Intro}
\emph{Robust optimization} approach considers the cases in which optimization problems often deal with uncertain data due to prediction errors, lack of information, fluctuations, or disturbances \cite{2,41,3,44,4,5}. Particularly, in most cases these problems depend on conflicting goals due to multiobjective decision makers which have different optimization criteria. So, the \emph{robust multiobjective optimization} is highly interesting in optimization theory and important in applications.

The first idea of robustness treated as a kind of sensitivity in the objective space against perturbations in the decision space for multiobjective optimization problems has been introduced in \cite{42,10}. Furthermore, some different concepts in minimax robustness for multiobjective optimization have been established in \cite{17,13,43}. More recently, the robustness concepts used in uncertain multiobjective optimization have been studied in \cite{23,38,18,39,40}.

To the best of our knowledge, the most powerful results in this direction were established for robust optimization in the finite-dimensional case. Hence, an infinite-dimensional framework would be suitable to investigate when involving optimality and duality in robust multiobjective optimization. From this, we are motivated to state and analyze problems that deal with infinite-dimensional frameworks.

Suppose that \mbox{$f: Z \to Y$} be a locally Lipschitzian vector-valued function between \emph{Asplund} spaces, and that $K \subset Y$ be a pointed (i.e., $K \bigcap \,(-K) = \{0\}$) closed convex cone. We consider the following \emph{multiobjective} optimization problem:
\begin{equation*}\hypertarget{P}{}
\begin{aligned}
  (\mathrm{P}) \qquad \min\nolimits_{K} \,\,\, &f(z) \\
                        \textrm{s.t.} \,\,\, &g_{i}(z) \le 0, \quad i = 1, 2, \dots, n, \nonumber \\
\end{aligned}
\end{equation*}
where the functions $g_{i}: Z \to \mathbb{R}$, $i = 1, 2, \dots, n$, define the constraints. Problem (\hyperlink{P}{P}) in the face of data \emph{uncertainty} in the constraints can be captured by the following \emph{uncertain multiobjective} optimization problem:
\begin{equation*}\hypertarget{UP}{}
\begin{aligned}
  (\mathrm{UP}) \qquad \min\nolimits_{K} \,\,\, &f(z) \\
                        \textrm{s.t.} \,\,\, &g_{i}(z, u) \le 0, \quad i = 1, 2, \dots, n, \nonumber \\
\end{aligned}
\end{equation*}
where $z \in Z$ is the vector of \emph{decision} variable, $u$ is the vector of \emph{uncertain} parameter and $u \in \mathcal{U}$ for some \emph{sequentially compact} topological space $\mathcal{U}$, and $g_{i}: Z \times \mathcal{U} \to \mathbb{R}$, $i = 1, 2, \dots, n$, are given functions.

One of the powerful deterministic structures to study problem (\hyperlink{UP}{UP}) is the \emph{robust} optimization, which is known as the problem that the uncertain objective and constraint are satisfied for all possible scenarios within a prescribed uncertainty set. We now associate with them:
\begin{equation*}\hypertarget{RP}{}
\begin{aligned}
  (\mathrm{RP}) \qquad \min\nolimits_{K} \,\,\, &f(z) \\
                        \textrm{s.t.} \,\,\, &g_{i}(z, u) \le 0, \quad \forall u \in \mathcal{U}, \,\, i = 1, 2, \dots, n. \nonumber \\
\end{aligned}
\end{equation*}
The feasible set $F$ of problem (\hyperlink{RP}{RP}) is defined by
\begin{equation*}
  F := \big\{ z \in Z \,\mid\, g_{i}(z, u) \le 0, \,\, \forall u \in \mathcal{U}, i = 1, 2, \dots, n \big\}.
\end{equation*}

The following concepts of solutions can be found in the literature; see e.g., \cite{13,17}.

\begin{defn}\hypertarget{Def2-1}{}
\begin{itemize}
  \item [(i)] We say that a vector $\bar{z} \in Z$ is a \emph{robust efficient solution} of problem (\hyperlink{UP}{UP}), denoted by $\bar{z} \in \mathcal{S}(RP)$, if $\bar{z}$ is a \emph{efficient solution} of problem (\hyperlink{RP}{RP}), i.e., $\bar{z} \in F$ and
      \begin{equation*}
        f(z) - f(\bar{z}) \notin - K \setminus \{0\}, \quad \forall z \in  F.
      \end{equation*}
  \item [(ii)] A vector $\bar{z} \in Z$ is called a \emph{weakly robust efficient solution} of problem (\hyperlink{UP}{UP}), denoted by $\bar{z} \in \mathcal{S}^{w}(RP)$, if $\bar{z}$ is a \emph{weakly efficient solution} of problem (\hyperlink{RP}{RP}), i.e., $\bar{z} \in F$ and
      \begin{equation*}
        f(z) - f(\bar{z}) \notin - \textrm{int}\hspace{.4mm}K, \quad \forall z \in  F.
      \end{equation*}
  \item [(iii)] A vector $\bar{z} \in Z$ is called a \emph{properly robust efficient solution} of problem (\hyperlink{UP}{UP}), denoted by $\bar{z} \in \mathcal{S}^{pr}(RP)$, if $\bar{z}$ is a \emph{properly efficient solution} of problem (\hyperlink{RP}{RP}), i.e., $\bar{z} \in F$ and there exists $y^{*} \in \textrm{int}\hspace{.4mm}K^{+}$ such that
      \begin{equation*}
        \langle y^{*}, f(z) - f(\bar{z}) \rangle \ge 0, \quad \forall z \in  F.
      \end{equation*}
  \end{itemize}
\end{defn}

The problem (\hyperlink{UP}{UP}) with the decision vector taken in a Asplund space and uncertain parameters which lies in a sequentially compact topological space covers a wide range of robust multiobjective optimization problems (cf. \cite{26,22,8,25,6}). Our main purpose in this paper is to investigate a nonsmooth/nonconvex multiobjective optimization problem with uncertain constraints in arbitrary Asplund spaces under the pseudo convexity assumptions. We first establish necessary optimality theorem for weakly robust efficient solutions of problem (\hyperlink{UP}{UP}) by employing a fuzzy optimality condition of a nonsmooth/nonconvex multiobjective optimization problem without any constrained qualification in the sense of the Fr\'{e}chet subdifferential, and then derive necessary optimality condition for properly robust efficient solutions of the problem. Sufficient conditions for weakly and properly efficient solutions as well as for robust efficient solutions to such a problem are also provided by means of applying the new concepts of generalized pseudo convex functions. Along with optimality conditions, we address a Mond-Weir-type robust dual problem to problem (\hyperlink{UP}{UP}) and explore weak, strong, and converse duality properties under assumptions of pseudo convexity.

The outline of the paper is organized as follows. Section \ref{Sec2-Preli} recalls some preliminary definitions and several auxiliary results. In Section \ref{Sec3-NecSuf}, necessary/sufficient optimality conditions for weakly and properly robust efficient solutions and also sufficient condition for robust efficient solutions of problem (\hyperlink{UP}{UP}) are established in terms of the limiting subdifferential. Section \ref{Sec4-Duality} is devoted to presenting duality relations between the corresponding problems.

\section{Preliminaries}\label{Sec2-Preli}
Our notation is basically standard in the area of variational analysis; see, e.g., \cite{27}. Throughout this paper, all the spaces under consideration are \emph{Asplund}, unless otherwise stated, with the norm $\|\cdot\|$ and the canonical pairing $\langle \cdot\,,\cdot \rangle$ between the space $Z$ in question and its \emph{dual} $Z^{*}$ equipped with the \emph{weak$^{*}$ topology $w^{*}$}. By $B_{Z}(z,r)$, we denote the \emph{closed ball} centered at $z \in Z$ with radius $r > 0$, while $B_{Z}$ and $B_{Z^{*}}$ stand for the \emph{closed unit ball} in $Z$ and $Z^{*}$, respectively. Given a nonempty set $ \Gamma \subset Z $, the symbols $\textrm{co}\hspace{.4mm}\Gamma$, $\textrm{cl}\hspace{.4mm}\Gamma$, and $\textrm{int}\hspace{.4mm}\Gamma$ indicate the \emph{convex hull}, \emph{topological closure}, and \emph{topological interior} of $\Gamma$, respectively, while $\textrm{cl}^{*}\Gamma$ stands for the \emph{weak$^{*}$ topological closure} of $\Gamma \subset  Z^{*}$. The \emph{dual cone} of $\Gamma$ is the set
\begin{equation*}
  \Gamma^{+} := \big\{ z^{*} \in Z^{*} \,\mid\, \langle z^{*}, z \rangle \ge 0, \,\,\, \forall z \in \Gamma \big\}.
\end{equation*}
Besides, $\mathbb{R}^{n}_{+}$ signifies the nonnegative orthant of $\mathbb{R}^{n}$ for $n \in \mathbb{N} := \{1,2,\dots\}$.

A given set-valued mapping $H : \Gamma \subset Z \overrightarrow{\to} Z^{*}$ is called \emph{weak$^{*}$ closed} at $\bar{z} \in \Gamma$ if for any sequence $\{z_{k}\} \subset \Gamma$, $z_{k} \to \bar{z}$, and any sequence $\{z^{*}_{k}\} \subset Z^{*}$, $z^{*}_{k} \in H(z_{k})$, $z^{*}_{k} \overset{w^{*}} \to z^{*}$, one has $z^{*} \in H(\bar{z})$.

For a set-valued mapping $H : Z \overrightarrow{\to} Z^{*}$, the \emph{sequential Painlev\'{e}-Kuratowski upper/outer limit} of $H$ as $z \to \bar{z}$ is defined by
\begin{align*}
  \underset{z \to \bar{z}}{\textrm{Lim}\sup} \, H(z) := \Big\{ z^{*} \in Z^{*} \,\mid\,\,\, &\exists
  \text{ sequences } z_{k} \to \bar{z} \text{ and } z^{*}_{k} \overset{{\scriptscriptstyle w^{*}}} \to z^{*} \\
  &\text{with } z^{*}_{k} \in H(z_{k}) \text{ for all } k \in \mathbb{N} \Big\}.
\end{align*}

Let $\Gamma \subset Z$ be \emph{locally closed} around $\bar{z} \in \Gamma$, i.e., there is a neighborhood $U$ of $\bar{z}$ for which $\Gamma \bigcap \textrm{cl}\hspace{.4mm}U$ is closed. The \emph{Fr\'{e}chet normal cone} $\widehat{N}(\bar{z}; \Gamma)$ and the \emph{Mordukhovich normal cone} $N(\bar{z}; \Gamma)$ to $\Gamma$ at $\bar{z} \in \Gamma$ are defined by
\begin{align}
  \label{2-1}
  \widehat{N}(\bar{z}; \Gamma) &:= \Big\{z^{*} \in Z^{*} \,\mid\, \limsup\limits_{z \overset{\hspace{-1mm}\sGamma} \to \bar{z}} \dfrac{\langle z^{*}, z - \bar{z} \rangle}{\|z - \bar{z}\|} \le 0\Big\}, \\
  \label{2-2}
  N(\bar{z}; \Gamma) &:= \underset{z \overset{\hspace{-1mm}\sGamma} \to \bar{z}}{\textrm{Lim}\sup} \, \widehat{N}(z; \Gamma),
\end{align}
where $z \overset{\hspace{-1mm}\Gamma} \to \bar{z}$ stands for $z \to \bar{z}$ with $z \in \Gamma$. If $\bar{z} \notin \Gamma$, we put $\widehat{N}(\bar{z}; \Gamma) = N(\bar{z}; \Gamma) := \emptyset$.

For an extended real-valued function $\psi : Z \to \overline{\mathbb{R}}$, the \emph{limiting/Mordukhovich subdifferential} and the \emph{regular/Fr\'{e}chet subdifferential} of $\psi$ at $\bar{z} \in \textrm{dom}\,\psi$ are given, respectively, by
\begin{equation*}
  \partial \psi(\bar{z}) := \big\{ z^{*} \in Z^{*} \,\mid\, (z^{*}, -1) \in  N((\bar{z}, \psi(z)); \textrm{epi}\,\psi) \big\}
\end{equation*}
and
\begin{equation*}
  \widehat{\partial} \psi(\bar{z}) := \big\{ z^{*} \in Z^{*} \,\mid\, (z^{*}, -1) \in  \widehat{N}((\bar{z}, \psi(z)); \textrm{epi}\,\psi) \big\}.
\end{equation*}
If $|\psi(\bar{z})| = \infty$, then one puts $\partial \psi(\bar{z}) := \widehat{\partial} \psi(\bar{z}) := \emptyset$.

Put $\langle y^{*}, f \rangle (z) := \langle y^{*}, f(z) \rangle$, $z \in Z$, $y^{*} \in Y^{*}$, for a vector-valued map $f : Z \to Y$, and denote $\textrm{gph}\,f := \big\{ (z, y) \in Z \times Y \,\mid\, y = f(z) \big\}$. Next we recall the needed results known as the scalarization formulae of the \emph{coderivatives}.

\begin{lem}\label{Lem2-1}
Let $y^{*} \in Y^{*}$, and let $f : Z \to Y$ be Lipschitz around $\bar{z} \in Z$. We have
\begin{itemize}
  \item [\emph{(i)}] \emph{(See \cite[Proposition~3.5]{28})} $z^{*} \in \widehat{\partial} \langle y^{*}, f \rangle(\bar{z}) \,\, \Leftrightarrow \,\, (z^{*}, -y^{*}) \in \widehat{N}((\bar{z}, f(\bar{z})); \text{gph}\,\,f)$.
  \item [\emph{(ii)}] \emph{(See \cite[Theorem~1.90]{27})} $z^{*} \in \partial \langle y^{*}, f \rangle(\bar{z}) \,\, \Leftrightarrow \,\, (z^{*}, -y^{*}) \in N((\bar{z}, f(\bar{z})); \text{gph}\,\,f)$.
\end{itemize}
\end{lem}

Another calculus result is the \emph{sum rule} for the limiting subdifferential.

\begin{lem}\label{Lem2-2} \emph{(See \cite[Theorem~3.36]{27})}
Let $ \psi_{i} : Z \to \overline{\mathbb{R}}$, $(i \in \{1, 2,\dots, n\}, n \ge 2)$, be lower semicontinuous around $\bar{z}$, and let all but one of these functions be Lipschitz continuous around $\bar{z} \in Z$. Then, one has
\begin{equation*}
  \partial ( \psi_{1} + \psi_{2} +\dots+ \psi_{n})(\bar{z}) \subset \partial \psi_{1}(\bar{z}) + \partial \psi_{2}(\bar{z}) + \dots + \partial \psi_{n}(\bar{z}).
\end{equation*}
\end{lem}

The following lemma computes the limiting subdifferential for the \emph{maximum} functions in Asplund spaces. The interested reader is referred to \cite{25,50,51} for more details and proofs. The notation $\partial_{z}$ signifies the limiting subdifferential operation with respect to $z$.

\begin{lem}\label{Lem2-6}
Let $\mathcal{U}$ be a sequentially compact topological space, and let $g : Z \times \mathcal{U} \to \mathbb{R}$ be a function such that for each fixed $u \in \mathcal{U}$, $g(\cdot, u)$ is locally Lipschitz on $U \subset Z$ and for each fixed $z \in U$, $g(z, \cdot)$ is upper semicontinuous on $\mathcal{U}$. Let $\psi(z) := \max\limits_{u \in \mathcal{U}} g(z, u)$. If the multifunction $(z, u) \in U \times \mathcal{U} \,\, \overrightarrow{\to} \,\, \partial_{z} g(z, u) \subset Z^{*}$ is weak$^{*}$ closed at $(\bar{z}, \bar{u})$ for each $\bar{u} \in \mathcal{U}(\bar{z})$, then the set \mbox{$\emph{cl}^{*}\emph{co} \Big(\medcup \Big\{\partial_{z} g(\bar{z}, u) \,\mid\, u \in \mathcal{U}(\bar{z})\Big\}\Big)$} is nonempty and
\begin{equation*}
\partial \psi(\bar{z}) \subset \emph{cl}^{*}\emph{co} \Big(\medcup \Big\{\partial_{z} g(\bar{z}, u) \,\mid\, u \in \mathcal{U}(\bar{z})\Big\}\Big),
\end{equation*}
where $\mathcal{U}(\bar{z}) = \big\{ u \in \mathcal{U} \,\mid\, g(\bar{z}, u) = \psi(\bar{z}) \big\} $.
\end{lem}

In what follows, we also use a formula for the limiting subdifferential of maximum of finitely many functions in Asplund spaces.
\begin{lem}\label{Lem2-7}\emph{(See \cite[Theorem~3.46]{27})} Let $\psi_{i} : Z \to \overline{\mathbb{R}}$, $(i \in \{1, 2,\dots, n\}, n \ge 2)$, be Lipschitz continuous around $\bar{z}$. Set $\psi(z) := \max\limits_{i \in \{1, 2,\dots, n\}} \psi_{i}(z)$. Then
\begin{equation*}
  \partial \psi(\bar{z}) \subset \medcup \Big\{\partial \Big(\sum_{i \in I(\bar{z})} \mu_{i} \, \psi_{i} \Big)(\bar{z}) \,\mid\, (\mu_{1}, \mu_{2},\dots, \mu_{n}) \in \Lambda(\bar{z})\Big\},
\end{equation*}
where
\begin{equation*}
  I(\bar{z}) := \big\{ i \in \{1, 2,\dots, n\} \,\mid\, \psi_{i}(\bar{z}) = \psi(\bar{z}) \big\}
\end{equation*}
and
\begin{equation*}
  \Lambda(\bar{z}) := \Big\{(\mu_{1}, \mu_{2},\dots, \mu_{n}) \,\mid\, \mu_{i} \ge 0, \,\, \sum_{i=1}^{n} \mu_{i} = 1, \,\, \mu_{i} \, (\psi_{i}(\bar{z}) - \psi(\bar{z})) = 0 \Big\}.
\end{equation*}
\end{lem}

Motivated by the concept of pseudo-quasi generalized convexity due to Fakhar \cite{8}, we introduce a similar
concept of pseudo-quasi convexity type for $f$ and $g$.

\begin{defn}\hypertarget{Def2-3}{}
  \begin{itemize}
    \item [(i)] $f$ is \emph{pseudo convex} at $\bar{z} \in Z$ if for any $z \in Z$ and $y^{*} \in K^{+}$ the following holds:
        \begin{equation*}
          \langle y^{*}, f \rangle(z) < \langle y^{*}, f \rangle(\bar{z}) \Longrightarrow  \big(\langle v^{*}, z - \bar{z} \rangle  < 0, \quad \forall v^{*} \in {\partial} \langle y^{*}, f \rangle(\bar{z})\big).
        \end{equation*}
    \item [(ii)] $f$ is \emph{strictly pseudo convex} at $\bar{z} \in Z$ if for any $z \in Z \setminus \{\bar{z}\}$ and $y^{*} \in K^{+} \setminus \{0\}$ the following holds:
        \begin{equation*}
          \langle y^{*}, f \rangle(z) \le \langle y^{*}, f \rangle(\bar{z}) \Longrightarrow  \big(\langle v^{*}, z - \bar{z} \rangle  < 0, \quad \forall v^{*} \in {\partial} \langle y^{*}, f \rangle(\bar{z})\big).
        \end{equation*}
    \item [(iii)] $g$ is \emph{generalized quasi convex} at $\bar{z} \in Z$ if for any $z \in Z$ and $u \in \mathcal{U}$ the following holds:
        \begin{equation*}
           g_{i}(z, u) \le g_{i}(\bar{z}, u) \Longrightarrow \big(\langle u^{*}_{i}, z - \bar{z} \rangle \le 0, \quad \forall u^{*}_{i} \in \partial_{z} g_{i}(\bar{z}, u) \big), \,\, i = 1, 2, \dots, n.
        \end{equation*}
  \end{itemize}
\end{defn}

\begin{defn}\hypertarget{Def2-4}{}
\begin{itemize}
\item [(i)] We say that $(f, g)$ is \emph{type I pseudo convex} at $\bar{z} \in Z$ if $f$ and $g$ are pseudo convex and generalized quasi convex at $\bar{z} \in Z$, respectively.
\item [(ii)] We say that $(f, g)$ is \emph{type II pseudo convex} at $\bar{z} \in Z$ if $f$ and $g$ are strictly pseudo convex and generalized quasi convex at $\bar{z} \in Z$, respectively.
\end{itemize}
\end{defn}

\begin{rem}\hypertarget{Exa2-1(Rem)}{}
\begin{itemize}
  \item [(i)] It follows from Definitions \hyperlink{Def2-3}{2.1} and \hyperlink{Def2-4}{2.2} that if $(f, g)$ is type II pseudo convex at $\bar{z} \in Z$, then $(f, g)$ is type I pseudo convex at $\bar{z} \in Z$, but converse is not true (see Example \ref{Exa2-2}).
  \item [(ii)] It is noted that the generalized (resp., strictly generalized) convexity (see \cite[Definition~3.9]{25}) of $(f, g)$ is reduced to the type I (resp., type II) pseudo convexity of $(f, g)$. Furthermore, as the next example demonstrates, the class of type I pseudo convex functions is properly wider than the class of generalized convex functions, which is properly wider than convex functions (see \cite[Example~3.10]{25}).
\end{itemize}
\end{rem}

\begin{exa}\label{Exa2-2}
Let $Z := \mathbb{R}^{2}$, $Y := \mathbb{R}^{3}$, $\mathcal{U} := [-1, -\dfrac{1}{4}]$, and let $K := \mathbb{R}^{3}_{+}$. Consider $f : Z \to Y$ and $g : Z \times \mathcal{U} \to \mathbb{R}^{2}$ defined by $f := (f_{1}, f_{2}, f_{3})$ and $g := (g_{1}, g_{2})$, respectively, where
\begin{equation*}
        \left\{\begin{aligned}
              f_{1}(z_{1},z_{2}) &= 5 |z_{1}| + \frac{2}{5} z_{2} + \frac{4}{5}, \\
              f_{2}(z_{1},z_{2}) &= \frac{1}{2} |z_{1}| + 3 z_{2} + 6, \\
              f_{3}(z_{1},z_{2}) &= 4 |z_{1}| + \frac{1}{2} z_{2} + 1
              \end{aligned}
        \right.
        \quad \text{and} \quad
        \left\{\begin{aligned}
              g_{1}(z_{1},z_{2}, u) &= \frac{1}{4} u^{2} |z_{1}| + \frac{1}{2} u^{2} z_{2} + \frac{1}{4} |u|, \\
              g_{2}(z_{1},z_{2}, u) &= \frac{1}{8} z_{1}^{2} + |u| z_{2} + |u| + \frac{1}{4},
              \end{aligned}
        \right.
\end{equation*}
$u \in \mathcal{U}$. Let $\bar{z} := (0, -2) \in Z$. It is easy to calculate directly by the definitions that
\begin{equation*}
  \partial f_{1}(\bar{z}) = [-5, 5] \times \{\dfrac{2}{5}\}, \,\,\,\, \partial f_{2}(\bar{z}) = [-\dfrac{1}{2}, \dfrac{1}{2}] \times \{3\}, \,\, \text{ and } \,\, \partial f_{3}(\bar{z}) = [-4, 4] \times \{\dfrac{1}{2}\}.
\end{equation*}
Moreover
\begin{equation*}
  \partial_{z} g_{1}(\bar{z}, u) = [-\dfrac{1}{4} u^{2}, \dfrac{1}{4} u^{2}] \times \{\dfrac{1}{2} u^{2}\} \,\, \text{ and } \,\, \partial_{z} g_{2}(\bar{z}, u) = (0, |u|)
\end{equation*}
for all $u \in \mathcal{U}$.

Suppose that for some $z := (z_{1}, z_{2}) \in Z$ and $y^{*} := (y^{*}_{1}, y^{*}_{2}, y^{*}_{3})^{\top} \in K$ the condition $\langle y^{*}, f \rangle(z) < \langle y^{*}, f \rangle(\bar{z})$ is satisfied. Hence
\begin{equation*}
  (\frac{2}{5} y^{*}_{1} + 3 y^{*}_{2} + \frac{1}{2} y^{*}_{3}) z_{2} < - (5 y^{*}_{1} + \frac{1}{2} y^{*}_{2} + 4 y^{*}_{3}) |z_{1}| - \frac{4}{5} y^{*}_{1} - 6 y^{*}_{2} - y^{*}_{3}.
\end{equation*}
Let us divide both sides of above inequality by $c := \dfrac{2}{5} y^{*}_{1} + 3 y^{*}_{2} + \dfrac{1}{2} y^{*}_{3} \ne 0$, and so
\begin{equation}\label{2-18}
  z_{2} < - \frac{1}{c} (5 y^{*}_{1} + \frac{1}{2} y^{*}_{2} + 4 y^{*}_{3}) |z_{1}| - 2.
\end{equation}
By using (\ref{2-18}), for any $v^{*} := y^{*}_{1} v^{*}_{1} + y^{*}_{2} v^{*}_{2} + y^{*}_{3} v^{*}_{3} \in \partial \langle y^{*}, f \rangle(\bar{z})$, where $v^{*}_{i} := (v^{*}_{iz}, v^{*}_{iy}) \in \partial f_{i}(\bar{z})$, $i = 1, 2, 3$, we have
\begin{align*}
  \langle v^{*}, z - \bar{z} \rangle & = \begin{pmatrix}
                                           y^{*}_{1} v^{*}_{1z} + y^{*}_{2} v^{*}_{2z} + y^{*}_{3} v^{*}_{3z} \\
                                           \dfrac{2}{5} y^{*}_{1} + 3 y^{*}_{2} + \dfrac{1}{2} y^{*}_{3}
                                         \end{pmatrix}
                                         \times
                                         \begin{pmatrix}
                                           z_{1} \\
                                           z_{2} + 2
                                         \end{pmatrix} \\
                                     & < (y^{*}_{1} v^{*}_{1z} + y^{*}_{2} v^{*}_{2z} + y^{*}_{3} v^{*}_{3z}) z_{1} - (5 y^{*}_{1} + \frac{1}{2} y^{*}_{2} + 4 y^{*}_{3}) |z_{1}| \\
                                     & \le 0,
\end{align*}
where the latter inequality is due to $v^{*}_{1z} \in [-5, 5]$, $v^{*}_{2z} \in [-\dfrac{1}{2}, \dfrac{1}{2}]$, and $v^{*}_{3z} \in [-4, 4]$. So $\langle v^{*}, z - \bar{z} \rangle < 0$.

Now, for $z \in Z$ and $u \in \mathcal{U}$ the condition $g_{1}(z, u) \le g_{1}(\bar{z}, u)$ implies $\dfrac{1}{2} u^{2} z_{2} \le -\dfrac{1}{4} u^{2} |z_{1}| - u^{2}$, and thus for any $u^{*}_{1} := (u^{*}_{1z}, u^{*}_{1y}) \in \partial_{z}g_{1}(\bar{z}, u)$, we get
\begin{equation*}
  \langle u^{*}_{1}, z - \bar{z} \rangle \le u^{*}_{1z} z_{1} - \frac{1}{4} u^{2} |z_{1}| \le 0
\end{equation*}
due to $u^{*}_{1z} \in [-\dfrac{1}{4} u^{2}, \dfrac{1}{4} u^{2}]$. Similarly for $z \in Z$ and $u \in \mathcal{U}$ satisfying $g_{2}(z, u) \le g_{2}(\bar{z}, u)$, one has $z_{2} \le - \dfrac{1}{8|u|} z_{1}^{2} - 2$, and thus for any $u^{*}_{2} \in \partial_{z}g_{2}(\bar{z}, u)$, it holds
\begin{equation*}
  \langle u^{*}_{2}, z - \bar{z} \rangle \le - \dfrac{1}{8} z_{1}^{2} \le 0.
\end{equation*}
Therefore, $(f, g)$ is type I pseudo convex at $\bar{z}$.

However, there exist $z := (1, -3) \in Z$ and $y^{*} := \Big(0, \dfrac{7}{5}, 1\Big)^{\top} \in K$ such that $\langle y^{*}, f \rangle (z) = 0 = \langle y^{*}, f \rangle (\bar{z})$, but for $v_{1}^{*} := \Big(0, \dfrac{2}{5}\Big)^{\top} \in \partial f_{1}(\bar{z})$, $v_{2}^{*} := \Big(\dfrac{1}{2}, 3\Big)^{\top} \in \partial f_{2}(\bar{z})$, and \mbox{$v_{3}^{*} := \Big(4, \dfrac{1}{2}\Big)^{\top} \in \partial f_{3}(\bar{z})$,} we have
\begin{equation*}
  v^{*} := y_{1}^{*} v_{1}^{*} + y_{2}^{*} v_{2}^{*} + y_{3}^{*} v_{3}^{*} = \Big(\dfrac{47}{10}, \dfrac{47}{10}\Big)^{\top}
\end{equation*}
so $\langle v^{*}, z - \bar{z} \rangle = 0$. This implies that $(f, g)$ is not type II pseudo convex at $\bar{z}$. On the other hand, there exist $z := (-1, -5) \in Z$ and $y^{*} := (0, 0, 0)^{\top} \in K$ such that for any $\nu \in Z$ one has
\begin{align*}
  \langle y^{*}, f \rangle (z) &- \langle y^{*}, f \rangle (\bar{z}) = 0, \\
  g_{1}(z, u) &- g_{1}(\bar{z}, u) = - \dfrac{5}{4} u^{2} < 0, \\
  g_{2}(z, u) &- g_{2}(\bar{z}, u) = \dfrac{1}{8} - 3 |u| < 0.
\end{align*}
Hence, $(f, g)$ is not type I convex at $\bar{z}$.
\end{exa}

\begin{exa}\label{Exa2-3}
Let $Z$, $Y$, $\mathcal{U}$, and $K$ be the same as Example \ref{Exa2-2}. Let $f : Z \to Y$ defined by $f := (f_{1}, f_{2}, f_{3})$, where
\begin{equation*}
        \left\{\begin{aligned}
              f_{1}(z_{1}, z_{2}) &= \frac{4}{5} z_{1}^{2} + 5|z_{1}| + \frac{4}{5} (z_{2} + 2)^{2} + \frac{2}{5} z_{2} + \frac{4}{5}, \\
              f_{2}(z_{1}, z_{2}) &= 6 z_{1}^{2} +\frac{1}{2} |z_{1}| + 6 (z_{2} + 2)^{2} + 3 z_{2} + 6, \\
              f_{3}(z_{1}, z_{2}) &= z_{1}^{2} + 4 |z_{1}| + (z_{2} + 2)^{2} + \frac{1}{2} z_{2} + 1,
              \end{aligned}
        \right.
\end{equation*}
and let $g : Z \times \mathcal{U} \to \mathbb{R}^{2}$ be the same as Example \ref{Exa2-2}. Let $\bar{z} := (0, -2) \in Z$. Then
\begin{equation*}
  \partial f_{1}(\bar{z}) = [-5, 5] \times \{\dfrac{2}{5}\}, \,\,\,\, \partial f_{2}(\bar{z}) = [-\dfrac{1}{2}, \dfrac{1}{2}] \times \{3\}, \,\, \text{ and } \,\, \partial f_{3}(\bar{z}) = [-4, 4] \times \{\dfrac{1}{2}\}.
\end{equation*}

Suppose that for some $z := (z_{1}, z_{2}) \in Z \setminus \{\bar{z}\}$ and $y^{*} := (y^{*}_{1}, y^{*}_{2}, y^{*}_{3})^{\top} \in K \setminus \{0\}$ the condition $\langle y^{*}, f \rangle(z) \le \langle y^{*}, f \rangle(\bar{z})$ is satisfied. Hence
\begin{align*}
  (\frac{2}{5} y^{*}_{1} + 3 y^{*}_{2} + \frac{1}{2} y^{*}_{3}) z_{2} \le & - (\frac{4}{5} y^{*}_{1} + 6 y^{*}_{2} + y^{*}_{3}) z_{1}^{2} - (5 y^{*}_{1} + \frac{1}{2} y^{*}_{2} + 4 y^{*}_{3}) |z_{1}| \\
  &- (\frac{4}{5} y^{*}_{1} + 6 y^{*}_{2} + y^{*}_{3}) (z_{2} + 2)^{2} - \frac{4}{5} y^{*}_{1} - 6 y^{*}_{2} - y^{*}_{3}.
\end{align*}
Let us divide both sides of above inequality by $c := \dfrac{2}{5} y^{*}_{1} + 3 y^{*}_{2} + \dfrac{1}{2} y^{*}_{3} \ne 0$, and so
\begin{equation}\label{2-22}
  z_{2} \le - 2 z_{1}^{2} - \dfrac{1}{c} (5 y^{*}_{1} + \dfrac{1}{2} y^{*}_{2} + 4 y^{*}_{3}) |z_{1}| - 2 (z_{2} + 2)^{2} - 2 .
\end{equation}
By using (\ref{2-22}), for any $v^{*} := y^{*}_{1} v^{*}_{1} + y^{*}_{2} v^{*}_{2} + y^{*}_{3} v^{*}_{3} \in \partial \langle y^{*}, f \rangle(\bar{z})$, where $v^{*}_{i} := (v^{*}_{iz}, v^{*}_{iy}) \in \partial f_{i}(\bar{z})$, $i = 1, 2, 3$, we have
\begin{align*}
  \langle v^{*}, z - \bar{z} \rangle & \le (y^{*}_{1} v^{*}_{1z} + y^{*}_{2} v^{*}_{2z} + y^{*}_{3} v^{*}_{3z}) z_{1} - 2 a z_{1}^{2} - (5 y^{*}_{1} + \frac{1}{2} y^{*}_{2} + 4 y^{*}_{3}) |z_{1}| - 2 a (z_{2} + 2)^{2} \\
 & < 0,
\end{align*}
where the latter strict inequality is due to $v^{*}_{1z} \in [-5, 5]$, $v^{*}_{2z} \in [-\dfrac{1}{2}, \dfrac{1}{2}]$, $v^{*}_{3z} \in [-4, 4]$, and $z \ne \bar{z}$. So $\langle v^{*}, z - \bar{z} \rangle < 0$. The complete calculations is similar to that of Example \ref{Exa2-2}. Therefore, $(f, g)$ is type II pseudo convex at $\bar{z}$.
\end{exa}

\begin{assum}\hypertarget{assum}{}(\hspace{-.05mm}See \cite[p.131]{25})
Let $\mathcal{U}$ be a sequentially compact topological space, and let $g: Z \times \mathcal{U} \to \mathbb{R}^{n}$ be a function satisfying the following hypotheses:
\begin{itemize}
  \item[(A1)] For a fixed $\bar{z} \in Z$, $g$ is locally Lipschitz in the first argument and uniformly on $\mathcal{U}$ in the second argument, i.e., there exist an open neighborhood $U$ of $\bar{z}$ and a positive constant $\ell$ such that $\|g(x, u) - g(y, u)\| \le \ell \|x - y\|$ for all $x, y \in U$ and $u \in \mathcal{U}$.
  \item[(A2)] For each $i = 1, 2, \dots, n$, the function $u \in \mathcal{U} \mapsto g_{i}(z, u) \in \mathbb{R}$ is upper semicontinuous for each $z \in U$.
  \item[(A3)] For each $i = 1, 2, \dots, n$, we define real-valued functions $\psi_{i}, \psi : Z \to \mathbb{R}$ via
      \begin{equation*}
         \psi_{i}(z) := \max_{u \in \mathcal{U}} g_{i}(z, u) \,\,\,\, \text{ and } \,\,\,\, \psi(z) := \max_{i \in \{1, 2, \dots, n\}} \psi_{i}(z),
      \end{equation*}
      and we observe that above assumptions imply that $\psi_{i}$ is well defined on $\mathcal{U}$. Besides, $\psi_{i}$ and $\psi$ follow readily that are locally Lipschitz on $U$, since each $g_{i}(\cdot, u)$ is (see \cite[(H1), p.131]{25} and \cite[p.290]{6}). Note that the feasible set $F$ can be equivalently characterized by:
      \begin{equation*}
        F = \big\{z \in Z \,\mid\, \psi_{i}(z) \le 0, \,\, i = 1, 2, \dots, n\big\} = \big\{z \in Z \,\mid\, \psi(z) \le 0\big\}.
      \end{equation*}
  \item[(A4)] For each $i = 1, 2, \dots, n$, the multifunction $(z, u) \in U \times \mathcal{U} \,\, \overrightarrow{\to} \,\, \partial_{z} g_{i}(z, u) \subset Z^{*}$ is weak$^{*}$ closed at $(\bar{z}, \bar{u})$ for each $\bar{u} \in \mathcal{U}_{i}(\bar{z})$, where $\mathcal{U}_{i}(\bar{z}) = \big\{ u \in \mathcal{U} \,\mid\, g_{i}(\bar{z}, u) = \psi_{i}(\bar{z}) \big\}$.
\end{itemize}
\end{assum}

In the rest of this section, we state a suitable constraint qualification in the sense of robustness, which is needed to obtain a so-called \emph{robust Karush-Kuhn-Tucker} (KKT) \emph{condition}.

\begin{defn}\hypertarget{Def2-5}{}(\hspace{-.05mm}See \cite[Definition 3.2]{25})
Let $\bar{z} \in F$. We say that \emph{constraint qualification} (CQ) \emph{condition} is satisfied at $\bar{z}$ if
\begin{equation*}
  0 \notin \textrm{cl}^{*}\textrm{co}\Big(\medcup \Big\{\partial_{z} g_{i}(\bar{z}, u) \,\mid\, u \in \mathcal{U}_{i}(\bar{z})\Big\}\Big),\quad i \in I(\bar{z}),
\end{equation*}
where $I(\bar{z}) := \big\{ i \in \{1, 2,\dots, n\} \,\mid\, \psi_{i}(\bar{z}) = \psi(\bar{z}) \big\}$.
\end{defn}

It is worth to mention here that this condition (CQ) reduces to the \emph{extended Mangasarian-Fromovitz constraint qualification} in the \emph{smooth} setting; see e.g., \cite{27} for more details.

\begin{defn}\hypertarget{Def2-6}{}
A point $\bar{z} \in F$ is said to satisfy the \emph{robust} (KKT) \emph{condition} if there exist $y^{*} \in K^{+} \setminus \{0\}$, $\mu := (\mu_{1}, \mu_{2},\dots, \mu_{n}) \in \mathbb{R}^{n}_{+}$, and $\bar{u}_{i} \in \mathcal{U}$, $i=1,2,\dots,n$, such that
\begin{equation*}
      \left\{\begin{aligned}
             & 0 \in \partial \langle y^{*}, f \rangle(\bar{z}) + \sum_{i=1}^{n} \mu_{i} \,  \textrm{cl}^{*}\textrm{co} \Big(\medcup \Big\{\partial_{z} g_{i}(\bar{z}, u) \,\mid\, u \in \mathcal{U}_{i}(\bar{z})\Big\}\Big), \\
             & \mu_{i} \, \max_{u \in \mathcal{U}} g_{i}(\bar{z}, u) = \mu_{i} \, g_{i}(\bar{z}, \bar{u}_{i}) = 0, \quad i=1,2,\dots,n.
             \end{aligned}
     \right.
\end{equation*}
\end{defn}

Therefore, the robust (KKT) condition defined above is guaranteed by the constraint qualification (CQ).

\section{Robust necessary and sufficient optimality conditions}\label{Sec3-NecSuf}
This section is devoted to study necessary optimality conditions for weakly and properly robust efficient solutions of problem (\hyperlink{UP}{UP}) by exploiting the nonsmooth version of Fermat's rule, the sum rule for the limiting subdifferential and the scalarization formulae of the coderivatives, and to discuss sufficient optimality conditions for such solutions as well as for robust efficient solutions by imposing the pseudo convexity assumptions.

The first theorem in this section establishes a necessary optimality condition in terms of the limiting subdifferential for weakly robust efficient solutions of problem (\hyperlink{UP}{UP}). To prove this theorem, we need a fuzzy necessary optimality condition in the sense of the Fr\'{e}chet subdifferential for weakly robust efficient solutions of problem (\hyperlink{UP}{UP}) as follows.

\begin{lem}\label{Thm3-1}\emph{(See \cite[Theorem~3.1]{34})}
Let $\bar{z} \in \mathcal{S}^{w}(RP)$. Then for each $k \in \mathbb{N}$, there exist $z^{1k} \in B_{Z}(\bar{z}, \frac{1}{k})$, \mbox{$z^{2k} \in B_{Z}(\bar{z}, \frac{1}{k})$,} $y^{*}_{k} \in K^{+}$ with $\| y^{*}_{k} \| = 1$, and $\alpha_{k} \in \mathbb{R}_{+}$ such that
\begin{align*}
   & 0 \in \widehat{\partial} \langle y^{*}_{k}, f \rangle (z^{1k}) + \alpha_{k} \widehat{\partial} \psi (z^{2k}) + \dfrac{1}{k} B_{Z^{*}}, \nonumber \\
   & | \alpha_{k} \psi(z^{2k}) | < \dfrac{1}{k}.
\end{align*}
\end{lem}

\begin{thm}\label{Thm3-2}
Suppose that $g_{i}$, $i = 1, 2, \dots, n$, satisfy \emph{\bf Assumptions} \hyperlink{assum}{\emph{(A1)-(A4)}}. If $\bar{z} \in \mathcal{S}^{w}(RP)$, then there exist $y^{*} \in K^{+}$, $\mu := (\mu_{1},\mu_{2},\dots,\mu_{n}) \in \mathbb{R}^{n}_{+}$, with $ \|y^{*}\| + \|\mu\| = 1$, and $\bar{u}_{i} \in \mathcal{U}$, $i = 1, 2, \dots, n$, such that
\begin{equation}\label{3-13}
 \left\{\begin{aligned}
      & 0 \in \partial \langle y^{*}, f \rangle(\bar{z}) + \sum_{i = 1}^{n} \mu_{i} \, \emph{cl}^{*}\emph{co} \Big(\medcup \Big\{\partial_{z} g_{i}(\bar{z}, u) \,\mid\, u \in \mathcal{U}_{i}(\bar{z})\Big\}\Big), \\
      & \mu_{i} \, \max\limits_{u \in \mathcal{U}} g_{i}(\bar{z}, u) = \mu_{i} \, g_{i}(\bar{z}, \bar{u}_{i}) = 0, \quad i = 1, 2, \dots, n.
       \end{aligned}
 \right.
\end{equation}
Furthermore, if the \emph{(CQ)} is satisfied at $\bar{z}$, then \emph{(\ref{3-13})} holds with $y^{*} \ne 0$.
\end{thm}

\begin{proof}
Let $\bar{z} \in \mathcal{S}^{w}(RP)$. Applying Lemma \ref{Thm3-1}, we find sequences $z^{1k} \to \bar{z}$, $z^{2k} \to \bar{z}$, $y^{*}_{k} \in K^{+}$ with $\|y^{*}_{k}\| = 1$, $\alpha_{k} \in \mathbb{R}_{+}$, $z^{*}_{1k} \in \widehat{\partial} \langle y^{*}_{k}, f \rangle(z^{1k})$, and $z^{*}_{2k} \in \alpha_{k} \, \widehat{\partial} \psi(z^{2k})$ satisfying
\begin{align}
  \label{3-14}
  &0 \in z^{*}_{1k} + z^{*}_{2k} + \dfrac{1}{k} B_{Z^{*}}, \\
  &\alpha_{k} \, \psi(z^{2k}) \to 0 \text{ as } k \to \infty. \nonumber
\end{align}
To proceed, we consider the following two possibilities for the sequence $\{\alpha_{k}\}$:

\textbf{Case 1:} If $\{\alpha_{k}\}$ is bounded, there is no loss of generality in assuming that \mbox{$\alpha_{k} \to \alpha \in \mathbb{R}_{+}$} as $k \to \infty$. Moreover, since the sequence $\{y^{*}_{k}\} \subset K^{+}$ is bounded, by using the weak$^{*}$ sequential compactness of bounded sets in duals to Asplund spaces we may assume without loss of generality that $y^{*}_{k} \overset{w^{*}} \rightarrow \bar{y}^{*} \in K^{+}$ with $\|\bar{y}^{*}\| = 1$ as $k \to \infty$. Let $\ell > 0$ be a Lipschitz modulus of $f$ around $\bar{z}$. One has the inequalities $\|z^{*}_{1k}\| \le \ell \, \|y^{*}_{k}\| \le \ell$ for all $k \in \mathbb{N}$ (see, \cite[Proposition~1.85]{27}). In this way, by taking a subsequences if necessary that $z^{*}_{1k} \overset{w^{*}} \rightarrow z^{*}_{1} \in Z^{*}$ as $k \to \infty$ and thus it stems from (\ref{3-14}) that $z^{*}_{2k} \overset{w^{*}} \rightarrow z^{*}_{2} := - z^{*}_{1}$ as $k \to \infty$. Applying the part (i) of Lemma \ref{Lem2-1} to the inclusion $z^{*}_{1k} \in \widehat{\partial} \langle y^{*}_{k}, f \rangle(z^{1k})$ gives us the relation
\begin{equation*}
  (z^{*}_{1k}, -y^{*}_{k}) \in \widehat{N}((z^{1k}, f(z^{1k})); \text{gph}\,\,f), \quad k \in \mathbb{N}.
\end{equation*}
Passing there to the limit as $k \to \infty$ and using the definitions of normal cones (\ref{2-1}) and (\ref{2-2}), we get $(z^{*}_{1}, -\bar{y}^{*}) \in N((\bar{z}, f(\bar{z})); \text{gph}\,\,f)$, which is equivalent to
\begin{equation}\label{3-16}
  z^{*}_{1} \in \partial \langle \bar{y}^{*}, f \rangle(\bar{z}),
\end{equation}
due to the part (ii) of Lemma \ref{Lem2-1}. Similarly, we obtain $z^{*}_{2} \in \alpha \, \partial \psi(\bar{z})$. The latter inclusions with (\ref{3-16}) imply that
\begin{equation}\label{3-20}
  0 \in \partial \langle \bar{y}^{*}, f \rangle(\bar{z}) + \alpha \, \partial \psi(\bar{z}),
\end{equation}
by taking into account that $z^{*}_{2} = - z^{*}_{1}$. Invoking now the formula for the limiting subdifferential of maximum functions in Lemma \ref{Lem2-7} gives us
\begin{equation}\label{3-21}
  \partial \psi(\bar{z}) \subset \medcup \Big\{\partial \Big(\sum_{i \in I(\bar{z})} \mu_{i} \, \psi_{i}\Big)(\bar{z}) \,\mid\, (\mu_{1},\mu_{2},\dots,\mu_{n})\in \Lambda(\bar{z})\Big\},
\end{equation}
where $I(\bar{z}) = \big\{ i \in \{1, 2, \dots, n\} \,\mid\, \psi_{i}(\bar{z}) = \psi(\bar{z}) \big\}$ and
\begin{equation*}
  \Lambda(\bar{z}) = \Big\{ (\mu_{1},\mu_{2},\dots,\mu_{n}) \,\mid\, \mu_{i} \ge 0, \,\, \sum_{i=1}^{n} \mu_{i} = 1, \,\, \mu_{i} \, (\psi_{i}(\bar{z}) - \psi(\bar{z})) = 0 \Big\}.
\end{equation*}
Employing further Lemma \ref{Lem2-6}, we arrive at
\begin{equation}\label{3-22}
  \partial \psi_{i}(\bar{z}) \subset \textrm{cl}^{*}\textrm{co} \Big(\medcup \Big\{\partial_{z} g_{i}(\bar{z}, u) \,\mid\, u \in \mathcal{U}_{i}(\bar{z})\Big\}\Big), \quad i = 1, 2, \dots, n,
\end{equation}
where $\mathcal{U}_{i}(\bar{z}) = \big\{ u \in \mathcal{U} \,\mid\, g_{i}(\bar{z}, u) = \psi_{i}(\bar{z}) \big\}$ and the set $ \textrm{cl}^{*}\textrm{co} \Big(\medcup \Big\{\partial_{z} g_{i}(\bar{z}, u) \,\mid\, u \in \mathcal{U}_{i}(\bar{z})\Big\}\Big)$ is nonempty. It follows from the sum rule of Lemma \ref{Lem2-2} and the relations (\ref{3-20})-(\ref{3-22}) that
\begin{equation*}
  0 \in \partial \langle \bar{y}^{*}, f \rangle(\bar{z}) +  \alpha \medcup \Big\{\sum_{i\in I(\bar{z})} \mu_{i} \, \textrm{cl}^{*}\textrm{co} \Big(\medcup \Big\{\partial_{z} g_{i}(\bar{z}, u) \,\mid\, u \in \mathcal{U}_{i}(\bar{z})\Big\}\Big) \,\mid\, (\mu_{1},\mu_{2},\dots,\mu_{n})
  \in \Lambda(\bar{z})\Big\}.
\end{equation*}
So, there exist $\bar{\mu} := (\bar{\mu}_{1},\bar{\mu}_{2},\dots,\bar{\mu}_{n}) \in \Lambda(\bar{z})$, with $\mathlarger{\sum\limits}_{i=1}^{n} \, \bar{\mu}_{i} = 1$ and $\bar{\mu}_{i} = 0$ for all $i \in \{1, 2, \dots, n\} \setminus I(\bar{z})$, such that
\begin{equation*}
 0 \in \partial \langle \bar{y}^{*}, f \rangle(\bar{z}) + \alpha \sum_{i=1}^{n} \bar{\mu}_{i} \, \textrm{cl}^{*}\textrm{co} \Big(\medcup \Big\{\partial_{z} g_{i}(\bar{z}, u) \,\mid\, u \in \mathcal{U}_{i}(\bar{z})\Big\}\Big).
\end{equation*}
Let us divide the above inclusion by $c := \|\bar{y}^{*}\| + \alpha \, \|\bar{\mu}\|$, and then put $y^{*} := \dfrac{\bar{y}^{*}}{c}$ and $\mu := \dfrac{\alpha}{c} \,\bar{\mu}$. Therefore, there exist $y^{*} \in K^{+}$ and $\mu := (\mu_{1}, \mu_{2},\dots, \mu_{n}) \in \mathbb{R}^{n}_{+}$, with $\|y^{*}\| + \|\mu\| = 1 $, such that
\begin{equation} \label{3-24}
 0 \in \partial \langle y^{*}, f \rangle(\bar{z}) + \sum_{i=1}^{n} \mu_{i} \, \textrm{cl}^{*}\textrm{co} \Big(\medcup \Big\{\partial_{z} g_{i}(\bar{z}, u) \,\mid\, u \in \mathcal{U}_{i}(\bar{z})\Big\}\Big).
\end{equation}

On the other side, by the sequentially compactness of $\mathcal{U}$ \mbox{and the upper semicontinuity of the} function $u \in \mathcal{U} \longmapsto g_{i}(\bar{z}, u)$ for each $i = 1, 2, \dots, n$, we can select $\bar{u}_{i} \in \mathcal{U}$ such that $g_{i}(\bar{z}, \bar{u}_{i}) = \max\limits_{u \in \mathcal{U}} g_{i}(\bar{z}, u) = \psi_{i}(\bar{z})$. Moreover, $\alpha \, \psi(\bar{z}) = 0$  due to $\alpha_{k} \, \psi(z^{2k}) \to 0$ as $k \to \infty$. Taking into account that $\psi_{i}(\bar{z}) = \psi(\bar{z})$ for all $i \in I(\bar{x})$, we obtain
\begin{equation*}
\mu_{i} \, g_{i}(\bar{z}, \bar{u}_{i}) = \dfrac{\alpha}{c} \, \bar{\mu}_{i} \, \psi_{i}(\bar{z}) =
\dfrac{\bar{\mu}_{i}}{c} \, [\alpha \, \psi(\bar{z})] = 0,
\end{equation*}
i.e., $\mu_{i} \, g_{i}(\bar{z}, \bar{u}_{i}) = \mu_{i} \, \max\limits_{u \in \mathcal{U}} g_{i}(\bar{z}, u) = 0$ for all $i \in \{1, 2, \dots, n\}$. This together with (\ref{3-24}) yields (\ref{3-13}).

\textbf{Case 2:} Assuming next that $\{\alpha_{k}\}$ is unbounded. Similar to the Case 1, we get from the inclusion $z^{*}_{2k} \in \alpha_{k} \, \widehat{\partial} \psi(z^{2k})$ that $(z^{*}_{2k}, - \alpha_{k}) \in \widehat{N}((z^{2k}, \psi(z^{2k})); \text{gph}\,\,\psi)$ for each $k \in \mathbb{N}$. So
\begin{equation*}
\Big(\dfrac{z^{*}_{2k}}{\alpha_{k}}, - 1\Big) \in \widehat{N}((z^{2k}, \psi(z^{2k})); \text{gph}\,\,\psi), \quad k \in \mathbb{N}.
\end{equation*}
Letting $k \to \infty$ and noticing (\ref{2-2}) again, we arrive at $(0, -1) \in N((\bar{z}, \psi(\bar{z})); \text{gph}\,\,\psi)$, which is equivalent to $0 \in \partial \psi(\bar{z})$. Proceeding as in the proof of the Case 1, we find $\mu := (\mu_{1},\mu_{2},\dots,\mu_{n}) \in \mathbb{R}^{n}_{+} \setminus \{0\}$, with $\|\mu\| = 1 $, satisfying
\begin{equation*}
 0 \in \sum_{i=1}^{n} \mu_{i} \, \textrm{cl}^{*}\textrm{co} \Big(\medcup \Big\{\partial_{z} g_{i}(\bar{z}, u) \,\mid\, u \in \mathcal{U}_{i}(\bar{z})\Big\} \Big).
\end{equation*}
We can also select $\bar{u}_{i} \in \mathcal{U}$ such that $\mu_{i} \, g_{i}(\bar{z}, \bar{u}_{i}) = \mu_{i} \, \psi_{i}(\bar{z}) =  \mu_{i} \, \psi(\bar{z}) = 0$ for each $i = 1, 2, \dots, n$, due to the unboundedness of $\{\alpha_{k}\}$ and $\alpha_{k} \, \psi(z^{2k}) \to 0$ as $k \to \infty$. So, (\ref{3-13}) holds by taking $y^{*} := 0 \in K^{+}$.

Finally, let $\bar{z}$ satisfy the (CQ) in the Case 1. It follows directly from (\ref{3-13}) that $y^{*} \ne 0$, which justifies the last statement of the theorem and completes the proof.
\end{proof}

\begin{rem}\hypertarget{Exa3-1(Rem)}{}
\begin{itemize}
    \item[(i)] Theorem \ref{Thm3-2} reduces to \cite[Theorem~3.3]{25} in the case of finite-dimensional optimization. Note further that our approach here, which involves the fuzzy necessary optimality condition in the sense of the Fr\'{e}chet subdifferential and the inclusion formula for the limiting subdifferential of maximum functions in the setting of Asplund spaces, is totally different from the one presented in the aforementioned paper.
   \item[(ii)] Observe that the results obtained in \cite[Theorem~3.3]{6}, \cite[Theorem~3.1]{26}, and \cite[Proposition~3.9]{8} are expressed for problems containing the convex or Q-convexlike objective functions and the convex constraint systems in terms of the convex or Clarke subdifferentials, but the one in Theorem \ref{Thm3-2} is established for nonconvex problems in the framework of Asplund spaces by applying the limiting subdifferential.
\end{itemize}
\end{rem}

Similarly, the following theorem presents a necessary optimality condition for properly robust efficient solutions of problem (\hyperlink{UP}{UP}).

\begin{thm}\label{Thm3-3}
Suppose that all conditions of \emph{Theorem \ref{Thm3-2}} are satisfied. If $\bar{z} \in \mathcal{S}^{pr}(RP)$, then there exist $y^{*} \in K^{+}$, $\mu := (\mu_{1},\mu_{2},\dots,\mu_{n}) \in \mathbb{R}^{n}_{+}$, with $ \|y^{*}\| + \|\mu\| = 1$, and $\bar{u}_{i} \in \mathcal{U}$, $i = 1, 2, \dots, n$, such that
\begin{equation}\label{3-27}
 \left\{\begin{aligned}
      & 0 \in \partial \langle y^{*}, f \rangle(\bar{z}) + \sum_{i=1}^{n} \mu_{i} \, \emph{cl}^{*}\emph{co} \Big(\medcup \Big\{\partial_{z} g_{i}(\bar{z}, u) \,\mid\, u \in \mathcal{U}_{i}(\bar{z})\Big\}\Big), \\
      & \mu_{i} \, \max_{u \in \mathcal{U}} g_{i}(\bar{z}, u) = \mu_{i} \, g_{i}(\bar{z}, \bar{u}_{i}) = 0, \quad i = 1, 2, \dots, n, \\
      & \langle y^{*}, f \rangle(\bar{z}) = \min\limits_{z \in F} \langle y^{*}, f \rangle(z).
       \end{aligned}
 \right.
\end{equation}
Furthermore, if the \emph{(CQ)} is satisfied at $\bar{z}$, then \emph{(\ref{3-27})} holds with $y^{*} \ne 0$.
\end{thm}

\begin{proof}
Let $\bar{z} \in \mathcal{S}^{pr}(RP)$. Then, there exists $\bar{y}^{*} \in \textrm{int}\hspace{.4mm}K^{+}$ such that
\begin{equation}\label{3-28}
   \langle \bar{y}^{*}, f(z) - f(\bar{z}) \rangle \ge 0, \quad \forall z \in  F.
\end{equation}
Let us define a real-valued function $\omega$ on $Z$ as follows:
\begin{equation*}
  \omega(z) := \max \big\{\langle \bar{y}^{*}, f(z) - f(\bar{z}) \rangle, \psi(z) \big\}, \quad \forall z \in Z.
\end{equation*}
It can be easily verified that $\omega(\bar{z}) = 0 \, \le \, \omega(z)$ on $Z$, which means that $\bar{z}$ is a minimizer for $\omega$. Employing the \emph{generalized Fermat's rule} (see \cite[Proposition~1.114]{27}), we have
\begin{equation}\label{3-29}
  0 \in \partial \omega(\bar{z}).
\end{equation}
From Lemma \ref{Lem2-7}, we deduce that
\begin{equation*}
  \partial \omega(\bar{z}) \subset \medcup \Big\{\partial \big(\alpha_{1} \, \langle \bar{y}^{*}, f(\cdot) - f(\bar{z}) \rangle + \alpha_{2} \, \psi(\cdot)\big)(\bar{z}) \,\mid \, \alpha_{1},\alpha_{2} \ge 0, \,\, \alpha_{1} + \alpha_{2} = 1, \,\, \alpha_{2} \, \psi(\bar{z}) = 0 \Big\}.
\end{equation*}
Combining this with (\ref{3-29}) and applying the sum rule of the limiting subdifferential give us $(\bar{\alpha}_{1}, \bar{\alpha}_{2}) \in \mathbb{R}^{2}_{+}$ with $\bar{\alpha}_{1} + \bar{\alpha}_{2} = 1$, such that $\bar{\alpha}_{2} \, \psi(\bar{z}) = 0$ and
\begin{equation*}
  0 \in \bar{\alpha}_{1} \, \partial \langle \bar{y}^{*}, f \rangle(\bar{z}) + \bar{\alpha}_{2} \, \partial \psi(\bar{z}).
\end{equation*}
Proceeding as in the proof of Theorem \ref{Thm3-2}, we find $\bar{\mu} := (\bar{\mu}_{1},\bar{\mu}_{2},\dots, \bar{\mu}_{n}) \in \Lambda(\bar{z})$, with $\mathlarger{\sum\limits}_{i=1}^{n} \, \bar{\mu}_{i} = 1$ and $\bar{\mu}_{i} = 0$ for each $i \in \{1, 2, \dots, n\} \setminus I(\bar{z})$, such that
\begin{equation*}
 0 \in \bar{\alpha}_{1} \, \partial \langle \bar{y}^{*}, f \rangle(\bar{z}) + \bar{\alpha}_{2} \, \sum_{i=1}^{n} \bar{\mu}_{i} \, \textrm{cl}^{*}\textrm{co} \Big(\medcup \Big\{\partial_{z} g_{i}(\bar{z}, u) \,\mid\, u \in \mathcal{U}_{i}(\bar{z})\Big\}\Big).
\end{equation*}
Dividing this by $c := \bar{\alpha}_{1}\,\|\bar{y}^{*}\| + \bar{\alpha}_{2}\,\|\bar{\mu}\|$ and  putting $y^{*} := \dfrac{\bar{\alpha}_{1}}{c}\,\bar{y}^{*}$ and $\mu := \dfrac{\bar{\alpha}_{2}}{c} \,\bar{\mu}$, we entail that $y^{*} \in K^{+}$ and $\mu \in \mathbb{R}^{n}_{+}$, with $\| y^{*}\| + \|\mu\| = 1 $, satisfying the first relation in the theorem. Furthermore, one can find $\bar{u}_{i} \in \mathcal{U}$ so that $\mu_{i} \, g_{i}(\bar{z}, \bar{u}_{i}) = 0$ due to $\bar{\alpha}_{2} \, \psi(\bar{z}) = 0$. To finish the proof of the theorem, we observe that
\begin{equation*}
  \langle y^{*}, f \rangle(z) - \langle y^{*}, f \rangle(\bar{z}) = \dfrac{\bar{\alpha}_{1}}{c} \, \langle \bar{y}^{*}, f(z) - f(\bar{z}) \rangle \overset{(\ref{3-28})}{\ge} 0, \quad \forall z \in F,
\end{equation*}
i.e. $\langle y^{*}, f \rangle(\bar{z}) = \min\limits_{z \in F} \langle y^{*}, f \rangle(z)$. Note that under the (CQ) condition, we have $y^{*} \ne 0$.
\end{proof}

We next revisit two examples to illustrate Theorem \ref{Thm3-2} and Theorem \ref{Thm3-3} for an uncertain multiobjective optimization problem.

\begin{figure}
    \centering
    \includegraphics[width=0.52\textwidth]{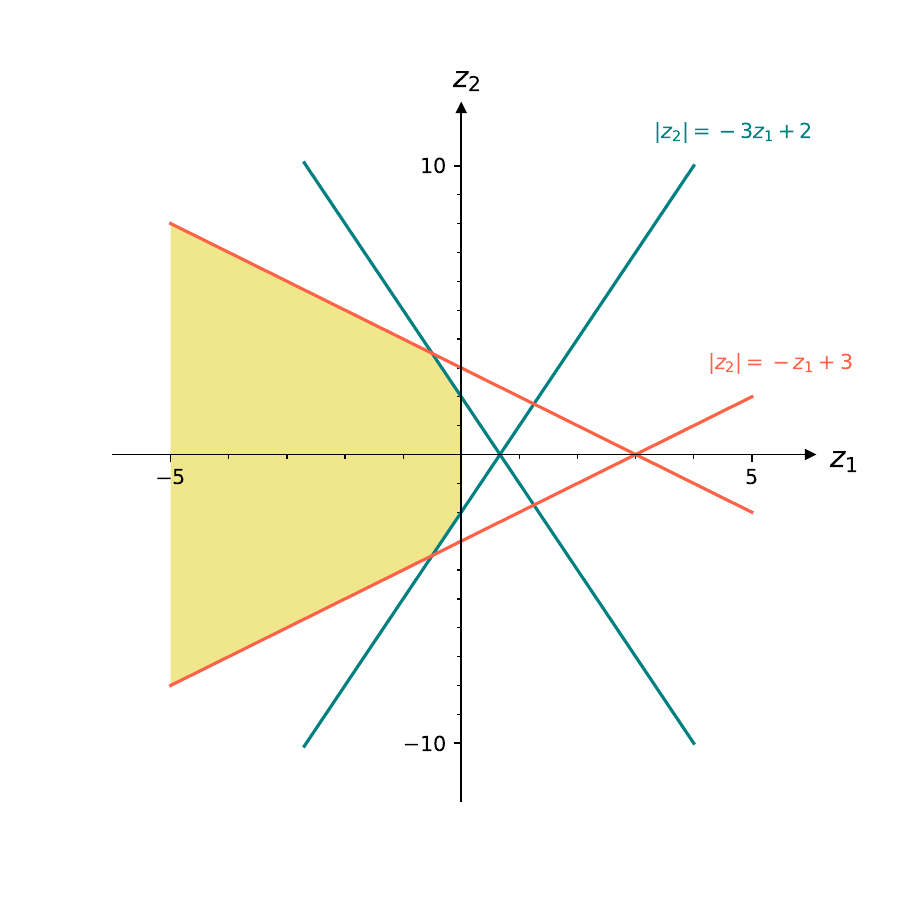}
    \vspace{-1cm}
    \caption{Robust feasible set of problem (UP) in Example \ref{Exa3-2}}
    \label{fig1}
\end{figure}

\begin{exa}\label{Exa3-2}
Suppose that $Z := \mathbb{R}^{2}$, $Y := \mathbb{R}^{3}$, $\mathcal{U} := [-1, 1]$, and $K := \mathbb{R}^{3}_{+}$. Consider the following uncertain optimization problem:
\begin{equation*}\hypertarget{UPexa3.2}{}
 (\mathrm{UP}) \qquad \min\nolimits_{K} \,\,\, \big\{ f(z) \,\mid\, g(z, u) \in - \mathbb{R}^{2}_{+} \big\},
\end{equation*}
where $f : Z \to Y$, $f := (f_{1}, f_{2}, f_{3})$ are defined by
\begin{equation*}
        \left\{\begin{aligned}
              f_{1} (z_{1},z_{2}) &= -2 z_{1} + |z_{2} - 1|,\\
              f_{2}(z_{1},z_{2}) &= \dfrac{1}{|z_{1}| + 1} - 3 z_{2} + 2,\\
              f_{3} (z_{1},z_{2}) &= \dfrac{1}{\sqrt{|z_{1}| + 1}} - |z_{2} - 1| - 1,
              \end{aligned}
        \right.
\end{equation*}
and $g : Z \times \mathcal{U} \to \mathbb{R}^{2}$, $g := (g_{1}, g_{2})$ are given by
\begin{equation*}
        \left\{\begin{aligned}
              g_{1}(z_{1},z_{2},u) &= u^{2} |z_{2}| + \max \big\{ z_{1}, 2 z_{1} \big\} - 3 |u|,\\
              g_{2}(z_{1},z_{2},u) &= -3 |z_{1}| + u z_{2} - 2,
              \end{aligned}
        \right.
\end{equation*}
where $u \in \mathcal{U}$. It is not difficult to check that
\begin{equation*}
  \big\{ u^{2} |z_{2}| + \max \big\{ z_{1}, 2 z_{1} \big\} - 3 |u| \le 0 \,\,\, \forall u \in \mathcal{U} \big\} = \{(z_{1},z_{2}) \in Z \,\mid\, z_{1} \le 0 \text{ and } |z_{2}| \le - z_{1} + 3\},
\end{equation*}
and, due to $z_{1} \le 0$, it can be verified that
\begin{equation*}
  \big\{ -3 |z_{1}| + u z_{2} - 2 \le 0 \,\,\, \forall u \in \mathcal{U} \big\} = \big\{ (z_{1},z_{2}) \in Z \,\mid\, z_{1} \le 0 \text{ and } |z_{2}| \le -3 z_{1} + 2 \big\}.
\end{equation*}
Therefore, the robust feasible set is
\begin{align*}
  F = &\big\{(z_{1}, z_{2}) \in Z \,\mid\, -\dfrac{1}{2} \le z_{1} \le 0 \text{ and } |z_{2}| \le -3 z_{1} + 2 \big\} \medcup \\
  &\big\{(z_{1}, z_{2}) \in Z \,\mid\, z_{1} \le -\dfrac{1}{2} \text{ and } |z_{2}| \le - z_{1} + 3 \big\},
\end{align*}
which is represented in Figure \ref{fig1}.

Let $\bar{z} := (0, 1) \in F$ and $z := (z_{1}, z_{2}) \in F$. Taking into account that $z_{1} \le 0$, we get $f_{1}(z) - f_{1}(\bar{z}) \ge 0$. So
\begin{equation*}
f(z) - f(\bar{z}) \notin - \text{int}\hspace{.4mm}K
\end{equation*}
for all $z \in F$, i.e., $\bar{z}$ is a weakly robust efficient solution of problem (\hyperlink{UPexa3.2}{UP}). Note further that
\begin{align*}
  \psi_{1}(\bar{z}) &= \max\limits_{u \in \mathcal{U}} g_{1}(\bar{z}, u) = \max\limits_{u \in \mathcal{U}} (u^{2} - 3 |u|) = 0, \\
  \psi_{2}(\bar{z}) &= \max\limits_{u \in \mathcal{U}} g_{2}(\bar{z}, u) = \max\limits_{u \in \mathcal{U}} (u - 2) = -1.
\end{align*}
So $\psi(\bar{z}) = \max \big\{ \psi_{1}(\bar{z}), \psi_{2}(\bar{z}) \big\} = 0$, $\mathcal{U}_{1}(\bar{z}) = \{0\}$, and $\mathcal{U}_{2}(\bar{z}) = \{1\}$. After calculations, we get
\begin{equation*}
\partial f_{1}(\bar{z}) = \{-2\} \times [-1, 1], \,\,\,\, \partial f_{2}(\bar{z}) = [-1, 1] \times \{-3\}, \,\,\,\, \partial f_{3}(\bar{z}) = \big[-\dfrac{1}{2}, \dfrac{1}{2}\big] \times \{-1, 1\},
\end{equation*}
and also
\begin{equation}\label{3-31}
        \left\{\begin{aligned}
              \textrm{cl}^{*}\textrm{co}\Big(\partial_{z}g_{1}(\bar{z}, u=0)\Big) &= [1, 2] \times \{0\},\\
              \textrm{cl}^{*}\textrm{co}\Big(\partial_{z}g_{2}(\bar{z}, u=1)\Big) &= [-3, 3] \times \{1\}.
              \end{aligned}
        \right.
\end{equation}
On the other side, since $I(\bar{z}) = \big\{ i \in \{1, 2\} \,\mid\, \psi_{i}(\bar{z}) = \psi(\bar{z}) \big\} = \{1\}$, it easily follows from (\ref{3-31}) that the (CQ) is satisfied at $\bar{z}$.

Finally, there exist $y^{*} = \Big(\dfrac{\sqrt{2}}{4}, 0, \dfrac{\sqrt{2}}{4}\Big)^{\top} \in K^{+}$ and $\mu = \Big(\dfrac{1}{2}, 0\Big)^{\top} \in \mathbb{R}^{2}_{+}$, with $\|y^{*}\| + \|\mu\| = 1$, such that
\begin{equation*}
  0 =
  \begin{pmatrix}
    \dfrac{\sqrt{2}}{4} & 0 & \dfrac{\sqrt{2}}{4}
  \end{pmatrix}
  \begin{pmatrix}
    -2 & 0 & 0 \\
    1 & -3 & -1
  \end{pmatrix}
  +
  \begin{pmatrix}
    \dfrac{1}{2} & 0
  \end{pmatrix}
  \begin{pmatrix}
    \sqrt{2} & 0 \\
    0 & 1
  \end{pmatrix},
\end{equation*}
and $\mu_{i} \, \max\limits_{u \in \mathcal{U}} g_{i}(\bar{z}, u) = 0$ for $i = 1, 2$.
\end{exa}

\begin{exa}\label{Exa3-3}
Let $Z$, $Y$, $\mathcal{U}$, and $K$ be the same as Example \ref{Exa3-2}. Consider the following uncertain optimization problem:
\begin{equation*}\hypertarget{UPexa3.3}{}
 (\mathrm{UP}) \qquad \min\nolimits_{K} \,\,\, \big\{ f(z) \,\mid\, g(z, u) \in - \mathbb{R}^{2}_{+} \big\},
\end{equation*}
where $f : Z \to Y$, $f = (f_{1}, f_{2}, f_{3})$ are defined by
\begin{equation*}
        \left\{\begin{aligned}
              f_{1}(z_{1},z_{2}) &= -2 z_{1} + |z_{2} - 1|,\\
              f_{2}(z_{1},z_{2}) &= \dfrac{1}{\sqrt{|z_{1}| + 1}} - 1,\\
              f_{3}(z_{1},z_{2}) &= -\dfrac{1}{\sqrt{2} (|z_{1}| + 1)} + |z_{2} - 1| + \dfrac{1}{\sqrt{2}}
              \end{aligned}
        \right.
\end{equation*}
and let $g : Z \times \mathcal{U} \to \mathbb{R}^{2}$ be the same as Example \ref{Exa3-2}.

Let $\bar{z} := (0, 1) \in F$. Observe that there exists $y^{*} = \Big(\dfrac{1}{5}, \dfrac{1}{5}, \dfrac{\sqrt{2}}{5} \Big)^{\top} \in \textrm{int}\hspace{.4mm}K^{+}$ such that
\begin{equation}\label{3-33}
  \langle y^{*}, f(z) - f(\bar{z}) \rangle = -\dfrac{2}{5} z_{1} + \Big(\dfrac{1+\sqrt{2}}{5}\Big) |z_{2} - 1| + \dfrac{1}{5} \Big(\dfrac{1}{\sqrt{|z_{1}| + 1}} - \dfrac{1}{|z_{1}| + 1}\Big) \overset{(z_{1} \le 0)} \ge 0
\end{equation}
for all $z \in F$, i.e., $\bar{z}$ is a properly robust efficient solution of problem (\hyperlink{UPexa3.3}{UP}). By the direct calculation, we have
\begin{equation*}
\partial f_{1}(\bar{z}) = \{-2\} \times [-1, 1], \,\,\,\, \partial f_{2}(\bar{z}) = [-\dfrac{1}{2}, \dfrac{1}{2}] \times \{0\}, \,\, \text{ and } \,\, \partial f_{3}(\bar{z}) = \{-\dfrac{\sqrt{2}}{2}, \dfrac{\sqrt{2}}{2}\} \times [-1, 1].
\end{equation*}
Again from Example \ref{Exa3-2}, the (CQ) holds at $\bar{z}$.

Now we find $\mu = \Big(\dfrac{3}{5},0\Big)^{\top} \in \mathbb{R}^{2}_{+}$, with $\|y^{*}\| + \|\mu\| = 1$, satisfying
\begin{equation*}
  0 =
  \begin{pmatrix}
    \dfrac{1}{5} & \dfrac{1}{5} & \dfrac{\sqrt{2}}{5}
  \end{pmatrix}
  \begin{pmatrix}
    -2 & 0 & -\dfrac{\sqrt{2}}{2} \\
    -\dfrac{\sqrt{2}}{2} & 0 & \dfrac{1}{2}
  \end{pmatrix}
  +
  \begin{pmatrix}
  \dfrac{3}{5} & 0
  \end{pmatrix}
  \begin{pmatrix}
    1 & 0 \\
    0 & 1
  \end{pmatrix},
\end{equation*}
and $\mu_{i} \, \max\limits_{u \in \mathcal{U}} g_{i}(\bar{z}, u) = 0$, $i = 1, 2$. Furthermore, (\ref{3-33}) implies $\min\limits_{z \in F} \langle y^{*}, f(z) \rangle = \langle y^{*}, f(\bar{z}) \rangle = 0$.
\end{exa}

The forthcoming theorem establishes a (KKT) sufficient optimality conditions for (weakly) robust efficient solutions of problem (\hyperlink{UP}{UP}).

\begin{thm}\label{Thm3-4}
Assume that $\bar{z} \in F$ satisfies the robust \emph{(KKT)} condition.
\begin{itemize}
  \item [\emph{(i)}] If $(f, g)$ is type I pseudo convex at $\bar{z}$, then $\bar{z} \in \mathcal{S}^{w}(RP)$.
  \item [\emph{(ii)}] If $(f, g)$ is type II pseudo convex at $\bar{z}$, then $\bar{z} \in \mathcal{S}(RP)$.
\end{itemize}
\end{thm}

\begin{proof}
Let $\bar{z} \in F$ satisfy the robust (KKT) condition. Therefore, there exist $y^{*} \in K^{+} \setminus \{0\}$, $v^{*} \in \partial \langle y^{*}, f \rangle(\bar{z})$, $\mu_{i} \ge 0$, and $u_{i}^{*} \in \textrm{cl}^{*}\textrm{co} \Big(\medcup \Big\{\partial_{z} g_{i}(\bar{z}, u) \,\mid\, u \in \mathcal{U}_{i}(\bar{z})\Big\}\Big)$, $i = 1, 2, \dots, n$, such that
\begin{align}\label{3-34}
  & 0 = v^{*} + \sum_{i=1}^{n} \mu_{i} \, u^{*}_{i}, \\
  \label{3-35}
  & \mu_{i} \, \max_{u \in \mathcal{U}} g_{i}(\bar{z}, u) = 0, \quad i = 1, 2, \dots, n.
\end{align}

Firstly, we justify (i). Argue by contradiction that $\bar{z} \notin \mathcal{S}^{w}(RP)$. Hence, there is $\hat{z} \in F$ such that $f(\hat{z}) - f(\bar{z}) \in - \textrm{int}\hspace{.4mm}K$. The latter gives us $\langle y^{*}, f(\hat{z}) - f(\bar{z}) \rangle < 0$ (see \cite[Lemma~3.21]{37}). Since $(f, g)$ is the type I pseudo convex at $\bar{z}$, we deduce from this inequality that
\begin{equation}\label{3-36}
  \langle v^{*}, \hat{z} - \bar{z} \rangle < 0.
\end{equation}
On the other side, it follows from (\ref{3-34}) for $\hat{z}$ above that
\begin{equation}\label{3-37}
  0 = \langle v^{*}, \hat{z} - \bar{z} \rangle + \sum_{i=1}^{n} \mu_{i} \, \langle u^{*}_{i}, \hat{z} - \bar{z} \rangle.
\end{equation}
The relations (\ref{3-36}) and (\ref{3-37}) entail that
\begin{equation*}
\sum_{i=1}^{n} \mu_{i} \, \langle u^{*}_{i}, \hat{z} - \bar{z} \rangle > 0.
\end{equation*}
To proceed, we assume that there is $i_{0} \in \{1, 2, \dots, n\}$ such that $\mu_{i_{0}} \, \langle u^{*}_{i_{0}}, \hat{z} - \bar{z} \rangle > 0$. Taking into account that $u^{*}_{i_{0}} \in \textrm{cl}^{*}\textrm{co} \Big(\medcup\Big\{\partial_{z} g_{i_{0}}(\bar{z}, u) \,\mid\, u \in \mathcal{U}_{i_{0}}(\bar{z})\Big\}\Big)$, we get sequence $\{u^{*}_{i_{0}k}\} \subset \textrm{co} \Big(\medcup \Big\{\partial_{z} g_{i_{0}}(\bar{z}, u) \,\mid\, u \in \mathcal{U}_{i_{0}}(\bar{z})\Big\}\Big)$ such that $u^{*}_{i_{0}k} \overset{{\scriptscriptstyle w^{*}}} \to u^{*}_{i_{0}}$. Hence, due to $\mu_{i_{0}} > 0$, there is $k_{0} \in \mathbb{N}$ such that
\begin{equation}\label{3-38}
  \langle u^{*}_{i_{0}k_{0}}, \hat{z} - \bar{z} \rangle > 0.
\end{equation}
In addition, since $u^{*}_{i_{0}k_{0}} \in \textrm{co} \Big(\medcup \Big\{\partial_{z} g_{i_{0}}(\bar{z}, u) \,\mid\, u \in \mathcal{U}_{i_{0}}(\bar{z})\Big\}\Big)$, there exist $u^{*}_{p} \in \medcup \Big\{\partial_{z} g_{i_{0}}(\bar{z}, u) \,\mid\, u \in \mathcal{U}_{i_{0}}(\bar{z})\Big\}$ and $\mu_{p} \ge 0$ with $\mathlarger{\sum\limits}_{p=1}^{s} \, \mu_{p} = 1$, $p = 1,2,\dots,s$, $s \in \mathbb{N}$, such that $u^{*}_{i_{0}k_{0}} = \mathlarger{\sum\limits}_{p=1}^{s} \, \mu_{p} \, u^{*}_{p}$. Combining the latter together (\ref{3-38}), we arrive at $\mathlarger{\sum\limits}_{p=1}^{s} \, \mu_{p} \, \langle u^{*}_{p}, \hat{z} - \bar{z} \rangle > 0$. Thus, we can take $p_{0} \in \{1,2,\dots,s\}$ such that
\begin{equation}\label{3-39}
  \langle u^{*}_{p_{0}}, \hat{z} - \bar{z} \rangle > 0,
\end{equation}
and choose $u_{i_{0}} \in \mathcal{U}_{i_{0}}(\bar{z})$ satisfying $u^{*}_{p_{0}} \in \partial_{z} g_{i_{0}}(\bar{z}, u_{i_{0}})$ due to $u^{*}_{p_{0}} \in \medcup \Big\{\partial_{z} g_{i_{0}}(\bar{z}, u) \,\mid\, u \in \mathcal{U}_{i_{0}}(\bar{z})\Big\}$. Invoking now definition of type I pseudo convexity of $(f, g)$ at $\bar{z}$, we get from (\ref{3-39}) that
\begin{equation}\label{3-40}
  g_{i_{0}}(\hat{z}, u_{i_{0}}) > g_{i_{0}}(\bar{z}, u_{i_{0}}).
\end{equation}
Note that $u_{i_{0}} \in \mathcal{U}_{i_{0}}(\bar{z})$, thus we have $g_{i_{0}}(\bar{z}, u_{i_{0}}) = \max\limits_{u \in \mathcal{U}} g_{i_{0}}(\bar{z}, u)$ which together with (\ref{3-35}) yields $\mu_{i_{0}} \, g_{i_{0}}(\bar{z}, u_{i_{0}}) = 0$. This implies by (\ref{3-40}) that $\mu_{i_{0}} \, g_{i_{0}}(\hat{z}, u_{i_{0}}) > 0$, and hence $g_{i_{0}}(\hat{z}, u_{i_{0}}) > 0$, which contradicts with the fact that $\hat{z} \in F$ and completes the proof of (i).

Assertion (ii) is proved similarly to the part (i). If $\bar{z} \notin \mathcal{S}(RP)$, then there exists $\hat{z} \in F$ such that $f(\hat{z}) - f(\bar{z}) \in - K \setminus \{0\}$. Therefore $\hat{z} \ne \bar{z}$ and $\langle y^{*}, f(\hat{z}) - f(\bar{z}) \rangle \le 0$. Now by using the definition of type II pseudo convexity of $(f, g)$  at $\bar{z}$, we arrive at the result.
\end{proof}

We immediately get the following sufficient optimality condition from Remark \hyperlink{Exa2-1(Rem)}{2.1}(i) and Theorem \ref{Thm3-4}.

\begin{cor}\label{Thm3-5(Cor)}
Let $\bar{z} \in F$ satisfy the robust \emph{(KKT)} condition and $(f, g)$ be type I pseudo convex at $\bar{z}$, then $\bar{z} \in \mathcal{S}(RP)$.
\end{cor}

\begin{rem}\hypertarget{Exa3-4(Rem)}{}
Theorem \ref{Thm3-4} improves \cite[Theorem~3.2]{5} and \cite[Theorem~3.11]{25}, where the involved functions are continuously differentiable and generalized convex in the setting of finite-dimensional spaces. We establish the (KKT) sufficient optimality condition for problem (\hyperlink{UP}{UP}) in the sense of pseudo convexity concept.
\end{rem}

Next, we prove a sufficient optimality theorem for the properly robust efficient solutions of problem (\hyperlink{UP}{UP}).

\begin{thm}\label{Thm3-6}
Let $(f, g)$ be type I pseudo convex at $\bar{z} \in F$. Assume that there exist \mbox{$y^{*} \in \emph{int}\hspace{.4mm}K^{+}$} and $\mu \in \mathbb{R}^{n}_{+}$ such that \emph{(\ref{3-13})} holds. Then $\bar{z} \in \mathcal{S}^{pr}(RP)$.
\end{thm}

\begin{proof}
Suppose that $\bar{z} \notin \mathcal{S}^{pr}(RP)$. There exists $\hat{z} \in F$ such that $\langle y^{*}, f(\hat{z}) - f(\bar{z}) \rangle < 0$. Following the definition of type I pseudo convexity of $(f, g)$ at $\bar{z}$, we get by the latter inequality that $\langle v^{*}, \hat{z} - \bar{z} \rangle < 0$. By using arguments similar to the proof of Theorem \ref{Thm3-4}(i), one can obtain the result.
\end{proof}

\section{Robust duality}\label{Sec4-Duality}
In this section, we formulate a Mond-Weir-type dual robust problem (\hyperlink{RD}{$\textrm{RD}_{MW}$}) for (\hyperlink{RP}{RP}), and investigate the weak, strong, and converse duality relations between the corresponding problems under pseudo convexity-affineness assumptions.

Let $y \in Z$, $y^{*} \in K^{+} \setminus \{0\}$, and $\mu \in \mathbb{R}^{n}_{+}$. In connection with the problem (\hyperlink{RP}{RP}), we introduce a \emph{dual robust multiobjective optimization} problem in the sense of Mond-Weir as follows:
\begin{equation*}\hypertarget{RD}{}
  (\textrm{RD}_{MW}) \qquad \max\nolimits_{K} \,\,\,\big\{ \bar{f}(y, y^{*}, \mu) := f(y) \,\mid\, (y, y^{*}, \mu) \in F_{MW} \big\}.
\end{equation*}
The feasible set $F_{MW}$ is given by
\begin{align*}
  F_{MW} &:= \bigg\{(y, y^{*}, \mu) \in Z \times K^{+} \setminus \{0\} \times \mathbb{R}^{n}_{+} \,\mid\, 0 \in \partial \langle y^{*}, f \rangle(y) + \sum_{i=1}^{n} \mu_{i} \, u^{*}_{i}, \\
  & u_{i}^{*} \in \textrm{cl}^{*}\textrm{co} \Big(\medcup \Big\{\partial_{z} g_{i}(y, u) \,\mid\, u \in \mathcal{U}_{i}(y) \Big\}\Big), \,\, \mu_{i} \, g_{i}(y, u) \ge 0, \,\, i = 1, 2, \dots, n \bigg\}.
\end{align*}

From now on, a robust efficient solution (resp., weakly robust efficient solution) of the dual problem (\hyperlink{RD}{$\textrm{RD}_{MW}$}) is defined similarly as in Definition \hyperlink{Def2-1}{1.1} by replacing $-K$ (resp., -\textrm{int}\hspace{.4mm}$K$) by $K$ (resp., \textrm{int}\hspace{.4mm}$K$). We denote the set of robust efficient solutions (resp., weakly robust efficient solutions) of problem (\hyperlink{RD}{$\textrm{RD}_{MW}$}) by $\mathcal{S}(RD_{MW})$ (resp., $\mathcal{S}^{w}(RD_{MW})$). Besides, we use the following notations for convenience:
\begin{align*}
  u \prec v \Leftrightarrow u-v \in -\textrm{int}\hspace{.4mm}K, \quad & u \nprec v \text{ is the negation of } u \prec v, \\
  u \preceq v \Leftrightarrow u-v \in -K\setminus\{0\}, \quad & u \npreceq v \text{ is the negation of } u \preceq v.
\end{align*}

The following theorem declares weak duality relations between the primal problem (\hyperlink{RP}{RP}) and the dual problem (\hyperlink{RD}{$\textrm{RD}_{MW}$}).

\begin{thm}\label{Thm4-1}\textbf{\textsc{(Weak Duality)}}
Let $z \in F$, and let $(y, y^{*}, \mu) \in F_{MW}$.
\begin{itemize}
  \item [\emph{(i)}] If $(f, g)$ is type I pseudo convex at $y$, then $f(z) \nprec \bar{f}(y, y^{*}, \mu)$.
  \item [\emph{(ii)}] If $(f, g)$ is type II pseudo convex at $y$, then $f(z) \npreceq \bar{f}(y, y^{*}, \mu)$.
  \end{itemize}
\end{thm}

\begin{proof}
By $(y, y^{*}, \mu) \in F_{MW}$, there exist $v^{*} \in \partial \langle y^{*}, f \rangle(y)$, $\mu_{i} \ge 0$, and $u_{i}^{*} \in \textrm{cl}^{*}\textrm{co} \Big(\medcup \Big\{\partial_{z} g_{i}(y, u) \,\mid\, u \in \mathcal{U}_{i}(y)\Big\}\Big)$, $i = 1, 2, \dots, n$, such that
\begin{align}\label{4-1}
  & 0 = v^{*} + \sum_{i=1}^{n} \mu_{i} \, u^{*}_{i},\\
  &  \mu_{i} \, g_{i}(y, u) \ge 0, \quad i = 1, 2, \dots, n.
\end{align}
To justify (i), assume that $f(z) \prec \bar{f}(y, y^{*}, \mu)$. Hence $\langle y^{*}, f(z) - \bar{f}(y, y^{*}, \mu) \rangle < 0$ due to $y^{*} \ne 0$. This is nothing else but $\langle y^{*}, f(z) - f(y) \rangle < 0$. Since $(f, g)$ is type I pseudo convex at $y$, we deduce from the last inequality that
\begin{equation*}
\langle v^{*}, z - y \rangle < 0.
\end{equation*}
On the other side, it follows from (\ref{4-1}) for $z$ above that
\begin{equation*}
  0 = \langle v^{*}, z - y \rangle + \sum_{i=1}^{n} \mu_{i} \, \langle u^{*}_{i}, z - y \rangle.
\end{equation*}
Combining the latter relations, we get that
\begin{equation*}
  \sum_{i=1}^{n} \mu_{i} \, \langle u^{*}_{i}, z - y \rangle > 0.
\end{equation*}
Now suppose that there is $i_{0} \in \{1, 2, \dots, n\}$ such that $\mu_{i_{0}} \, \langle u^{*}_{i_{0}}, z - y \rangle > 0$. Proceeding similarly to the proof of Theorem \ref{Thm3-4}(i) and replacing $\hat{z} - \bar{z}$ by $z - y$ give us $g_{i_{0}}(z, u_{i_{0}}) > 0$, which contradicts with $z \in F$.

Next to prove (ii), we proceed similarly to the part (i) by using the type II pseudo convexity of $(f, g)$ at $y$, if $f(z) \preceq \bar{f}(y, y^{*}, \mu)$, then $z \ne y$ and we deduce $\langle v^{*}, z - y \rangle < 0$.
\end{proof}

We now establish a strong duality theorem which holds between (\hyperlink{RP}{RP}) and (\hyperlink{RD}{$\textrm{RD}_{MW}$}).

\begin{thm}\label{Thm4-2}\textbf{\textsc{(Strong Duality)}}
Let $\bar{z} \in \mathcal{S}^{w}(RP)$ be such that the \emph{(CQ)} is satisfied at this point. Then, there exist $(\bar{y}^{*}, \bar{\mu}) \in K^{+} \setminus \{0\} \times \mathbb{R}^{n}_{+}$ such that $(\bar{z}, \bar{y}^{*}, \bar{\mu}) \in F_{MW}$. Furthermore,
\begin{itemize}
  \item [\emph{(i)}] If $(f, g)$ is type I pseudo convex at any $y \in Z$, then $(\bar{z}, \bar{y}^{*}, \bar{\mu}) \in \mathcal{S}^{w}(RD_{MW})$.
  \item [\emph{(ii)}] If $(f, g)$ is type II pseudo convex at any $y \in Z$, then $(\bar{z}, \bar{y}^{*}, \bar{\mu}) \in \mathcal{S}(RD_{MW})$.
\end{itemize}
\end{thm}

\begin{proof}
Thanks to Theorem \ref{Thm3-2}, we find $y^{*} \in K^{+} \setminus \{0\}$, $\mu_{i} \ge 0$, and $u_{i}^{*} \in \textrm{cl}^{*}\textrm{co} \Big(\medcup \Big\{\partial_{z} g_{i}(\bar{z}, u) \,\mid\, u \in \mathcal{U}_{i}(\bar{z}) \Big\}\Big)$, $i = 1, 2, \dots, n$, satisfying
\begin{align}
  & 0 \in \partial \langle y^{*}, f \rangle(\bar{z}) + \sum_{i=1}^{n} \mu_{i} \, u^{*}_{i}, \nonumber \\
  \label{4-6}
  & \mu_{i} \, \max_{u \in \mathcal{U}} g_{i}(\bar{z}, u) = 0, \quad i = 1, 2, \dots, n.
\end{align}
Putting $\bar{y}^{*} := y^{*}$ and $\bar{\mu} := (\mu_{1},\mu_{2},\dots,\mu_{n})$, we have $(\bar{y}^{*}, \bar{\mu}) \in K^{+} \setminus \{0\} \times \mathbb{R}^{n}_{+}$. Moreover, the inclusion $u \in \mathcal{U}_{i}(\bar{z})$ means that $g_{i}(\bar{z}, u) = \max\limits_{u \in \mathcal{U}} g_{i}(\bar{z}, u)$ for all $i \in \{1, 2, \dots, n\}$. Thus, it stems from (\ref{4-6}) that $\mu_{i} \, g_{i}(\bar{z}, u) = 0$, $i = 1, 2, \dots, n$. So $(\bar{z}, \bar{y}^{*}, \bar{\mu}) \in F_{MW}$.

(i) As $(f, g)$ be type I pseudo convex at any $y \in Z$, employing (i) of \mbox{Theorem \ref{Thm4-1} gives us}
\begin{equation*}
  \bar{f}(\bar{z}, \bar{y}^{*}, \bar{\mu}) = f(\bar{z}) \nprec \bar{f}(y, y^{*}, \mu)
\end{equation*}
for each $(y, y^{*}, \mu) \in F_{MW}$. Hence $(\bar{z}, \bar{y}^{*}, \bar{\mu}) \in \mathcal{S}^{w}(RD_{MW})$.

(ii) As $(f, g)$ be type II pseudo convex at any $y \in Z$, employing (ii) of Theorem \ref{Thm4-1} allows us
\begin{equation*}
  \bar{f}(\bar{z}, \bar{y}^{*}, \bar{\mu}) \npreceq \bar{f}(y, y^{*}, \mu)
\end{equation*}
for each $(y, y^{*}, \mu) \in F_{MW}$. Therefore $(\bar{z}, \bar{y}^{*}, \bar{\mu}) \in \mathcal{S}(RD_{MW})$.
\end{proof}

\begin{thm}\label{Thm4-3}\textbf{\textsc{(Strong Duality)}}
Let $\bar{z} \in \mathcal{S}^{pr}(RP)$ be such that the \emph{(CQ)} is satisfied at this point. Then, there exist $(\bar{y}^{*}, \bar{\mu}) \in K^{+} \setminus \{0\} \times \mathbb{R}^{n}_{+}$ such that $(\bar{z}, \bar{y}^{*}, \bar{\mu}) \in F_{MW}$. Furthermore,
\begin{itemize}
  \item [\emph{(i)}] If $(f, g)$ is type I pseudo convex at any $y \in Z$, then $(\bar{z}, \bar{y}^{*}, \bar{\mu}) \in \mathcal{S}^{w}(RD_{MW})$.
  \item [\emph{(ii)}] If $(f, g)$ is type II pseudo convex at any $y \in Z$, then $(\bar{z}, \bar{y}^{*}, \bar{\mu}) \in \mathcal{S}(RD_{MW})$.
\end{itemize}
\end{thm}

\begin{proof}
It follows from Theorem \ref{Thm3-3} and the proof of Theorem \ref{Thm4-2}.
\end{proof}

We finish this section by presenting converse duality relations between (\hyperlink{RP}{RP}) and (\hyperlink{RD}{$\textrm{RD}_{MW}$}).

\begin{thm}\label{Thm4-4}\textbf{\textsc{(Converse Duality)}}
Let $(\bar{z}, \bar{y}^{*}, \bar{\mu}) \in F_{MW}$ be such that $\bar{z} \in F$.
\begin{itemize}
   \item [\emph{(i)}] If $(f, g)$ is type I pseudo convex at $\bar{z}$, then $\bar{z} \in \mathcal{S}^{w}(RP)$.
   \item [\emph{(ii)}] If $(f, g)$ is type II pseudo convex at $\bar{z}$, then $\bar{z} \in \mathcal{S}(RP)$.
\end{itemize}
\end{thm}

\begin{proof}
Since $(\bar{z}, \bar{y}^{*}, \bar{\mu}) \in F_{MW}$, there exist $v^{*} \in \partial \langle \bar{y}^{*}, f \rangle(\bar{z})$, $\bar{\mu}_{i} \ge 0$, and $u_{i}^{*} \in \textrm{cl}^{*}\textrm{co} \Big(\medcup \Big\{\partial_{z} g_{i}(\bar{z}, u) \,\mid\, u \in \mathcal{U}_{i}(\bar{z})\Big\}\Big)$, $i = 1, 2, \dots, n$, such that
\begin{align}\label{4-7}
  & 0 = v^{*} + \sum_{i=1}^{n} \bar{\mu}_{i} \, u^{*}_{i}, \\
  & \bar{\mu}_{i} \, g_{i}(\bar{z}, u) \ge 0, \quad i = 1, 2, \dots, n. \nonumber
\end{align}
Let us justify (i). On the contrary, suppose that $\bar{z} \notin \mathcal{S}^{w}(RP)$. Hence, there is $\hat{z} \in F$ such that $f(\hat{z}) - f(\bar{z}) \in - \textrm{int}\hspace{.4mm}K$. The latter inclusion provides $\langle \bar{y}^{*}, f(\hat{z}) - f(\bar{z}) \rangle < 0$. By the type I pseudo convexity of $(f, g)$ at $\bar{z}$, we deduce from this inequality that
\begin{equation*}
\langle v^{*}, \hat{z} - \bar{z} \rangle < 0.
\end{equation*}
Furthermore, it follows from (\ref{4-7}) for $\hat{z}$ that
\begin{equation*}
  0 = \langle v^{*}, \hat{z} - \bar{z} \rangle + \sum_{i=1}^{n} \mu_{i} \, \langle u^{*}_{i}, \hat{z} - \bar{z} \rangle.
\end{equation*}
So, the above relationships entail that
\begin{equation*}
  \sum_{i=1}^{n} \mu_{i} \, \langle u^{*}_{i}, \hat{z} - \bar{z} \rangle > 0.
\end{equation*}
Now arguing as in the proof of Theorem \ref{Thm3-4}(i), one can arrive at the result.

The proof of (ii) is similar to that of (i), so we omit the corresponding details.
\end{proof}

\printbibliography
\small

\end{document}